\documentclass[11pt]{article}
\usepackage[english]{babel}
\usepackage[cp850]{inputenc}
\usepackage[dvips]{graphicx}
\usepackage{amsmath,amsfonts,amsthm,amssymb}
\usepackage[usenames, dvipsnames]{color}
\usepackage{fancyhdr}
\usepackage{stmaryrd}
\usepackage[colorlinks=true,citecolor=red,linkcolor=blue,urlcolor=blue,pdfstartview=FitH]{hyperref}
\usepackage{dsfont}
\usepackage{xcolor}
\usepackage{epsfig}

\font\tenmath=msbm10 scaled 1200

\font\sevenmath=msbm7 scaled 1200
\font\fivemath=msbm5 scaled 1200

\newfam\mathfam \textfont\mathfam=\tenmath
\scriptfont\mathfam=\sevenmath \scriptscriptfont\mathfam=\fivemath

\def\R{{\mathbb R}}
\def\N{{\mathbb N}}
\def\E{{\mathbb E}}

\def\P{{\mathbb P}}

\def\Q{{\mathbb Q}}

\def\D{{\mathbb D}}
\def\F{{\cal F}}

\newtheorem{Thm}{Theorem}[section]
\newtheorem{Lem}{Lemma}[section]
\newtheorem{Pro}{Proposition}[section]
\newtheorem{Cor}{Corollary}[section]
\newtheorem{Dfn}{Definition}[section]

\newfam\mathfam \textfont\mathfam=\tenmath
\scriptfont\mathfam=\sevenmath \scriptscriptfont\mathfam=\fivemath

\def \^#1{\if#1i{\accent"5E\i}\else{\accent"5E#1}\fi}

\def \D {I\!\!D}

\def \cqfd{\quad\Box}

\def \ss{\smallskip}
\def \bs{\bigskip}
\def \ni{\noindent}

\def\ni{\noindent}

\title{\bf Functional co-monotony of processes with an application to peacocks and barrier options}

\author{\textsc{Gilles Pag\`es} \thanks{Laboratoire de Probabilit\'es et Mod\`eles al\'eatoires, UMR~7599, UPMC, case 188, 4, pl. Jussieu, F-75252 Paris Cedex 5, France. E-mail: \texttt{gilles.pages@upmc.fr}}}

\date{September 4, 2012}

\begin{document}

\maketitle



\begin{abstract}We show that several general classes of stochastic processes satisfy a functional co-monotony principle, including processes with independent increments, Brownian bridge, Brownian diffusions, Liouville processes, fractional Brownian motion. As a first application, we recover some recent results about peacock processes obtained by Hirsch et al. in~\cite{YORetal} (see also~\cite{BOG}) which were themselves  motivated by a former work of  Carr et al. in~\cite{Carretal} about the sensitivities of Asian options with respect to their volatility and residual maturity (seniority). We also derive semi-universal bounds for various barrier options.
\end{abstract}



\section{Introduction}\label{Into}
The aim of this paper is to show that the classical co-monotony principle for real-valued random variables also holds for  large classes of stochastic processes like Brownian diffusion processes, processes with independent increments, Liouville processes, fractional Brownian motion(s), etc, if one considers the natural partial order on the space of real-valued functions defined on an interval. We also provide several examples of application, with a special emphasis on peacocks
(English quasi-acronym for ``processus croissants pour l'ordre convexe") inspired by recent works by Hirsch et al. in~\cite{YORetal}, which find themselves their original motivation in~\cite{Carretal} by Carr et al. about the sensitivities of Asian options in a Black-Scholes model. We also derive (semi-)universal upper or lower bound for various barrier options when the dynamics of the underlying asset price satisfies an appropriate functional  co-monotony  principle.

\smallskip
The starting point of what can be called {\em co-monotony principle} finds its origin in the following classical proposition dealing with one-dimensional real-valued random variables.

\begin{Pro}[One dimensional co-monotony principle] Let $X:(\Omega,{\cal A},\P)\to \R$ be a random variable and let $f,\, g:\R\to \R$ be two monotone  functions sharing the same monotony property.

\smallskip
\ni $(a)$ If $f(X)$, $g(X)$, $f(X)g(X)\!\in L^1(\P)$, then ${\rm Cov}(f(X),g(X))\ge 0$  $i.e.$
\[
\E \, f(X)g(X) \ge \E\, f(X)\E\, g(X).
\]
Furthermore, the inequality holds as an equality  if and only if $f(X)$ or $g(X)$ is $\P$-$a.s.$ constant.

\smallskip If   $f$ and $g$ are monotone with opposite monotony then the reverse inequality holds.

\smallskip
\ni  $(b)$ If $f$ and $g$ have  the same constant sign, then  integrability is no longer requested. As a consequence, if $f$ and $g$ have opposite monotony,  then
\[
 \|f(X)g(X)\|_1= \E\, f(X)g(X)\le \E\, f(X)\E \, g(X) =\|f(X)\|_1\|g(X)\|_1.
\]
\end{Pro}

These inequalities  are  straightforward consequences of Fubini's Theorem applied on $(\R\times\R, {\cal B}or(\R)^{\otimes 2}, \P_X^{\otimes 2})$ to the function $(x,x')\mapsto \big(f(x)-f(x')\big)\big(g(x)-g(x')\big)$ where $\P_X$ denotes the distribution of $X$.

\medskip
Typical applications of this scalar co-monotony principle are, among others, the antithetic simulation method  for variance reduction and more recently {\em a priori} sign results for the sensitivity of derivatives in Finance.

\medskip
\ni $\rhd$ {\em Antithetic simulation.}  Let $X:(\Omega,{\cal A},\P)\to \R$ be a random variable and let $\varphi:\R\to \R$ be a non-increasing function such that $\varphi(X)\stackrel{{\cal L}}{\sim} X$. Then, for every monotone function $f:\R\to \R$ such that $f(X)\!\in L^2(\P)$ and $\P\big(f(X)\neq \E f(X)\big)>0$, we have
\begin{eqnarray*}
{\rm Var}\Big(\frac{f(X)+f\circ \varphi(X)}{2}\Big) &=& \frac{2\big({\rm Var}(f(X))+{\rm Cov}(f(X), f\!\circ\!\varphi(X))\big)}{4} \\
&<& \frac{{\rm Var}(f(X))}{2}
\end{eqnarray*}
since ${\rm Cov}(f(X), f\!\circ\!\varphi(X))<0$. 

The variance is reduced by more than  a $2$-factor whereas the complexity of the simulation of  $\frac{f(X)+f\circ \varphi(X)}{2}$ is only twice higher than that of $f(X)$ (if one neglects the additional  cost  of the computation of $\varphi(x)$ compared to that of $x$).

\medskip
\ni $\rhd$ {\em Sensitivity (vega of an option).} Let $\varphi:(0,\infty)\to\R$ be a convex function with (at most) polynomial growth at $0$ and $+\infty$ in the sense that there exists a real constant $C>0$ such that 
$$
\forall\, x\!\in (0,+\infty),\qquad |\varphi(x)|\le C(x^r+x^{-r})
$$
and let $Z:(\Omega,{\cal A},\P)\to \R$ be an ${\cal N}(0;1)$-distributed random variable. Set for every $\sigma>0$
\[
f(\sigma) = \E\,\varphi\big(e^{\sigma Z-\frac{\sigma^2}{2}}\big).
\]
Although it does not appear as a straightforward consequence of its definition,  one easily derives from the above proposition 
that $f$ is a non-decreasing function of $\sigma$ on $(0,\infty)$. In fact, $\varphi$ is differentiable outside an at most countable subset of $(0,+\infty)$ (where its right and left derivatives differ) and its  derivative $\varphi'$ is non-decreasing, with polynomial growth as well since 
$$
|\varphi'(x)|\le \max\big(|\varphi(x+1)-\varphi(x)|, 2x^{-1}|\varphi(x)-\varphi(x/2)|\big),\;x\!\in (0,+\infty).
$$
  Since $Z$ has no atom, one easily checks that one can interchange derivative and expectation  to establish that $f$ is differentiable with derivative 
\[
f'(\sigma) = \E\Big(\varphi'\big(e^{\sigma Z-\frac{\sigma^2}{2}}\big)e^{\sigma Z-\frac{\sigma^2}{2}}(Z-\sigma)\Big),\quad \sigma>0.
\] 
A Cameron-Martin change of variable then yields 
\[
f'(\sigma)=\E\Big(\varphi'\big(e^{\sigma Z+\frac{\sigma^2}{2}}\big)Z\Big)
\]
so that, applying the co-monotony principle to the two non-decreasing (square integrable) functions  $z\mapsto \varphi'\big(e^{\sigma z-\frac{\sigma^2}{2}}\big)$ and $z\mapsto z$,  implies 
\[
f'(\sigma)\ge  \E\Big(\varphi'\big(e^{\sigma Z+\frac{\sigma^2}{2}}\big)\Big) \E\big( Z\big)=  \E\Big(\varphi'\big(e^{\sigma Z+\frac{\sigma^2}{2}}\big)\Big) \times 0 =0.
\]

\medskip
 Extensions of the above  co-monotony principle to functions on $\R^d$, $d\ge 2$, are almost as classical as the one dimensional case.  They can be established by induction when both functions $\Phi(x_1,\ldots,x_d)$ and  $\Psi(x_1,\ldots,x_d)$ defined on $\R^d$ are co-monotone in each variable $x_i$ ($i.e.$ having  the same or an opposite monotony property not depending on $i$) and when the $\R^d$-valued random vector $X$ has independent marginals.

Our aim in this paper is to show that this co-monotony principle can be again extended into a {\em  functional co-monotony principle}    satisfied by various classes of stochastic processes $X=(X_t)_{t\in [0,T]}$ whose paths lie in a sub-space $E$ of the vector space ${\cal F}([0,T],\R)$ of real-valued  functions defined on the interval $[0,T]$, $T>0$, equipped with the pointwise ({\em partial})  order on functions, defined by 
\[
\forall\, \alpha,\, \beta \in {\cal F}([0,T],\R),\; \alpha\le \beta \; \mbox{ if }\; \forall\, t\!\in [0,T],\; \alpha(t)\le \beta(t).
\]
Then a functional $F: E\to \R$ is said to be non-decreasing if
\[
\forall\, \alpha,\,\beta \!\in E, \quad \alpha\le \beta \Longrightarrow F(\alpha)\le F(\beta).
\]

The choice of $E$ will be motivated by the pathwise regularity of the process $X$. The space $E$ will also be endowed with a  metric topology  (and its Borel $\sigma$-field) so that $X$ can be seen as  an $E$-valued random vector. The functionals $F$ and $G$ involved in the co-monotony principle will be assumed to be continuous on $E$ (at least $\P_X$-$a.s.$).  Typical choices for $E$ will be $E={\cal C}([0,T], \R)$, ${\cal C}([0,T], \R^d)$, $I\!\!D ([0,T], \R)$ or $I\!\!D ([0,T], \R^d)$ and occasionally $L_{\R^d}^p([0,T], dt)$ (in this case we will switch to   the $dt$-$a.e.$ pointwise order instead of the pointwise order). Then by co-monotony principle, we mean that for every non-decreasing functionals  $F$ and $G$  defined on $E$, $\P_{X}$-$a.s.$ continuous,  
\[
\E\,F(X)G(X) \ge \E\, F(X)\E\, G(X).
\]
(The case of non-increasing functionals follows by considering the opposite functionals and the opposite monotony case by considering the opposite of only one of the functionals). Among the (classes of) processes of interest, we will consider continuous Gaussian processes with nonnegative covariance function (like the standard and the fractional Brownian motion, ``nonnegative"  Liouville processes), the  Markov processes with monotony preserving transitions (which includes of course Brownian diffusions), processes with independent increments, etc.

\smallskip The main problem comes from the fact that the naive pointwise order on functional spaces is not total so that the formal one-dimensional proof based on Fubini's theorem no longer applies.

\medskip As applications of such functional results, we will be able to extend the above sign property  for the {\em vega} of a ``vanilla"  option (whose payoff function is a function of the risky asset $S_T$ at the maturity $T$)  to ``exotic"   options. By ``exotic", we classically  mean that their  payoff 
is typically  a path-dependent   functional $F\big((S_t)_{t\in [0,T]}\big)$ of the risky asset   $(S_t)_{t\in [0,T]}$. The dynamics of this risky asset is still  a Black-Scholes model  where  $S^{\sigma}_t =s_0e^{\sigma W_t+(r-\frac{\sigma^2}{2})t}$, $s_0,\,\sigma\! >\! 0$. Doing so we will retrieve Carr et al. results about the sensitivity of Asian type options in a Black-Scholes model with respect to the volatility (see~\cite{Carretal}). We will also emphasize the close connection between co-monotony and the theory of {\em peacock}s (\footnote{stands for the French acronym PCOC (Processus Croissant pour l'Ordre Convexe).})   characterized by Kellerer~in~\cite{KEL} and recently put back into light in the book~\cite{YORetal} (see also the references therein). Let us briefly recall that an integrable  process $(X_{\lambda})_{\lambda \ge 0}$ is a peacock if   for  every convex function $\varphi:\R\to \R$, the $(-\infty,+\infty]$-valued 
function $\lambda\mapsto \E\, \varphi(X_{\lambda})$ is non-decreasing.  Kellerer's characterization theorem says that a process is a peacock if and only if there exists a martingale $(M_{\lambda})_{\lambda\ge0}$ such that $X_{\lambda}\stackrel{{\cal L}}{\sim} M_{\lambda}$, $\lambda\!\in \R_+$ (where $\stackrel{{\cal L}}{\sim} $ denotes  equality in distribution). Moreover, the process $(M_{\lambda})_{\lambda \ge 0}$ can be chosen to be Markovian. This proof being non-constructive, it does not help at all  establishing whether or not a process is a peacock. See also a new proof of Kellerer's Theorem due to  Hirsch and Roynette in~\cite{HIRO}. By contrast, one can find in~\cite{YORetal} a huge number of peacocks  with an explicit marginal martingale representation characterized through  various   tools from the theory of stochastic processes.

 More generally, when applied in its ``opposite" version, the co-monotony principle between nonnegative function simply provides a significant  improvement of  the   H\"older inequality since it makes   the $L^1$-norm sub-multiplicative. It can be used to produce less conservative bounds in various fields of applied probability, like recently in~\cite{LALEPA, ZBRA} where to provide bounds depending on functionals of a Brownian diffusion process, in the spirit of  the inequalities proposed in Section~\ref{bounds} for barrier options.
%

\medskip The paper is organized as follows:  Section~\ref{Deux} is devoted to the finite-dimensional co-monotony principle, Section~\ref{Trois}   to the functional co-monotony principle for continuous processes, Section~\ref{Quatre} to c\`adl\`ag processes. Section~\ref{Cinq} deals with  examples of applications,  to peacocks and to  exotic options for which we establish universal bounds (among price dynamics  sharing the functional co-monotony principle). 

\bigskip 
\ni{\sc Notation:} $\bullet$ $x_{0:n}=(x_0,\ldots,x_n)\!\in \R^{n+1}$, $x_{1:n}=(x_1,\ldots,x_n)\!\in \R^{n}$, etc. $(x|y)=\sum_{0\le k\le n} x_ky_k$ denotes the canonical inner product on $\R^{n+1}$. 

\smallskip \ni $\bullet$ We denote by $\le$ the componentwise order on $\R^{n+1}$ defined by $x_{0:n}\le x'_{0:n}$ if  $x_i\le x'_i$, $i=0,\ldots,n$. 

\smallskip \ni $\bullet$ ${\cal M}(d,r)$ denotes the vector space of matrices with $d$ rows and $r$ columns. $M^*$ denotes the transpose of matrix $M$.

\smallskip \ni  $\bullet$ $\perp\!\!\!\perp$ will emphasize in formulas the independence between two processes.
 
\smallskip \ni  $\bullet$ $\|\alpha\|_{\rm sup} = \sup_{t\in [0,T]}|\alpha(t)|$ for any function $\alpha:[0,T]\to \R$. 

\smallskip \ni  $\bullet$ $u_+$ denotes the positive part of the real number $u$. $\lambda_d$ denotes the Lebesgue measure on $(\R^d, {\cal B}or(\R^d))$ where ${\cal B}or(\R^d)$ denotes the Borel $\sigma$-field on $\R^d$.

\smallskip \ni $\bullet$ $X\stackrel{\cal L}{\sim} \mu$ means that  the random vector $X$ has distribution $\mu$.
\section{Finite-dimensional co-monotony principle}\label{Deux}

\subsection{Definitions and main results}

Let $(P(x,dy))_{x\in \R}$ be a  probability transition, $i.e.$ a family of probability measures such that for every $x\!\in \R$, $P(x,dy)$ is  a probability measure on  $(\R,{\cal B}or(\R))$ and  for every  $B\!\in{\cal B}or(\R)$, $x\mapsto P(x,B)$ is a Borel function.

\begin{Dfn}\label{Def:1} $(a)$   The transition $(P(x,dy))_{x\in \R}$ is monotony preserving if, for every bounded or nonnegative monotone function $f:\R\to\R$, the function $Pf$ defined for every real number $x\!\in \R$ by $Pf(x) = \int f(y)P(x,dy)$ is monotone with the same monotony.

\smallskip
\ni $(b)$ Two Borel  functions $\Phi$, $\Psi: \R^{d}\to \R$ are componentwise co-mono\-tone if, for every $i\!\in \{1,\ldots,d\}$ and every $d-1$-tuple $(x_1,\ldots,x_{i-1},x_{i+1},\ldots,x_d) \!\in \R^{d-1}$, both section functions $x_i\mapsto \Phi(x_1,\ldots,x_i,\ldots,x_d)$ and $x_i\mapsto \Psi(x_1,\ldots, x_i,\ldots,x_d)$ have the same monotony not depending on $(x_1,\ldots,x_{i-1},x_{i+1},\ldots,x_d)$  nor on $i$.

\smallskip
\ni $(c)$ If $\Phi$ and $-\Psi$ are co-monotone, $\Phi$ and $\Psi$ are said to be anti-monotone.
\end{Dfn}

\noindent {\bf Remark.} If $P$ is monotony preserving and $f:\R\to\R$ is a monotone function such that $f\!\in \cap_{x\in \R}L^1(P(x,dy))$, then $Pf$ has the same monotony  as $f$. This is an easy consequence of the Lebesgue dominated convergence theorem and the approximation of $f$ by the the ``truncated" bounded  functions $ f_N= (-N)\vee (f\wedge N)$ which have the same monotony as $f$.

\begin{Dfn}[Componentwise co-monotony principle]\label{Def:2} An $\R^d$-valued random vector $X$ satisfies a componentwise co-monotony principle if, for every pair of   Borel componentwise co-monotone functions $\Phi,\,\Psi: \R^d\to \R$ such that $\Phi(X)$, $\Psi(X)$, $\Phi(X)\Psi(X)\!\in L^1(\P)$, 
\begin{equation}\label{eq:comondimfin}
\E\, \Phi(X)\Psi(X) \ge \E \,\Phi(X)
\E \,\Psi(X).
\end{equation}
\end{Dfn}

\noindent {\bf Remarks.} $\bullet$ If $X$ satisfies a componentwise co-monotony principle, then, for every pair of componentwise anti-monotone functions $\Phi, \Psi:\R^{n+1}\to \R$, the reverse inequality holds. In both cases, if  $\Phi$ and $\Psi$ both take values in $\R_+$ or $\R_-$ then  the inequalities remain true (in $\overline \R$) without integrability assumption. 

\smallskip
\noindent $\bullet$ Owing  to elementary approximation arguments, it is clear that it suffices to check~(\ref{eq:comondimfin}) for  bounded or nonnegative componentwise co-monotone functions.

\bigskip As a  straightforward consequence of the fact that the  functions $x_{0:n}\mapsto x_k$ and $x_{0:n}\mapsto x_\ell$ are  co-monotone, we derive the following   necessary condition for the componentwise co-monotony property.

\begin{Pro} Let $X=(X_k)_{1\le k\le d}$ be a a random vector. If $X_k\!\in L^2$, $k=1,\ldots,d$, then,
\[
\forall\, k,\; \ell\!\in \{1,\ldots,d\}, \quad {\rm Cov}(X_k,X_{\ell})= \E\, X_kX_{\ell}- \E\, X_k\E\,X_{\ell}\ge 0
\]
$i.e.$ the covariance matrix of $X$ has nonnegative entries. 
\end{Pro}

In finite dimension, the main result on componentwise co-monotony is the following. 
\begin{Pro}\label{ProFinidim} $(a)$ Let $X=(X_k)_{0\le k\le n}$ be an $\R$-valued  Markov chain defined on a probability space $(\Omega,{\cal A},\P)$ having a  (regular) version of its transitions 
\[
P_{k-1,k}(x,dy) = \P(X_k\!\in dy\,|\, X_{k-1}=x),\; k=1,\ldots,n
\]
which are  monotony preserving in the above sense. Then $X$ satisfies a componentwise co-monotony principle. 

\smallskip
\ni $(b)$ If the random variables $X_0,\ldots,X_n$ are independent, the conclusion remains true under the following  {\em weak} co-monotony assumption: there exists a permutation $\tau$ of the index set $\{0,\ldots,n\}$ such that for every $i\!\in \{0,\ldots, n\}$ and every $(x_{0},\ldots,x_{i-1},x_{i+1},\ldots,x_n) \!\in \R^{n}$, $x_i\mapsto \Phi(x_0,\ldots, x_i,\ldots, x_n)$ and $x_i\mapsto\Psi(x_0,\ldots, x_i,\ldots,x_n)$ have the same monotony, possibly depending on $(x_{\tau(0)},\ldots, x_{\tau(i-1)})$. Then the same conclusion as in $(a)$ holds true.
\end{Pro}

\ni{\bf Proof.} $(a)$ One proceeds by induction on $n\!\in \N$. If $n=0$, the result follows from the scalar co-monotony principle applied to $X_0$ (with distribution $\mu_0$).

\smallskip \ni $(n)\Longrightarrow (n+1)$: We may assume that $\Phi$ and $\Psi$ are bounded and, by changing if necessary  the functionals into their opposite,  that they are both  componentwise non-decreasing. Put ${\cal F}^X_k=\sigma(X_0,\ldots,X_k)$, $k=0,\ldots,n$. It follows from the Markov property that 
\[
\E\big(\Phi(X_{0:n+1})\,|\,\F^X_n\big) = \Phi^{(n)}(X_{0:n}) 
\]
where 
\[
\Phi^{(n)}(x_{0:n})=P_{n,n+1}\big(\Phi(x_{0:n},.)\big)(x_n).
\]
In particular, we have $\E\big(\Phi(X_{0:n+1})\big)= \E \big(\Phi^{(n)}(X_{0:n})\big)$. 
Let  $x_{0:n}\!\in \R^{n+1}$. Applying the one dimensional co-monotony principle with the probability distribution $P_{n,n+1}(x_n,dy)$ to $\Phi(x_{0:n},.)$ and $\Psi(x_{0:n},.)$ we get  
\begin{eqnarray}
\nonumber (\Phi\,\Psi)^{(n)}(x_{0:n}) &=&  \;P_{n,n+1}\big(\Phi\,\Psi(x_{0:n},.)\big)(x_n)\\
\nonumber  &\ge& P_{n,n+1}\big(\Phi(x_{0:n},.)\big)(x_n)P_{n,n+1}\big(\Psi(x_{0:n},.)\big)(x_n)\\
\label{MarkovfiniD}
&=& \Phi^{(n)}(x_{0:n})\Psi^{(n)}(x_{0:n})
\end{eqnarray}
so that, considering $X_{0:n+1}$ and  taking expectation, we get
\begin{eqnarray*}
\E\, (\Phi\,\Psi)(X_{0:n+1}) &=& \E\big((\Phi\,\Psi)^{(n)}(X_{0,n})\big) \\
&\ge& \E\big( \Phi^{(n)}(X_{0:n})\Psi^{(n)}(X_{0:n})\big).
\end{eqnarray*}
  It is clear that, for every $i\!\in \{0,\ldots,n-1\}$, $x_i\mapsto \Phi^{(n)}(x_0,\ldots, x_n)$ is non-decreasing since the transition $P_{n,n+1}$  is a nonnegative operator. Now let $x_n,\,x'_n \!\in \R$, $x_n\le x'_n$. Then
\begin{eqnarray*}
P_{n,n+1}\big(\Phi(x_0,\ldots,x_n,.)\big)(x_n) & \le & P_{n,n+1}\big(\Phi(x_0,\ldots,x'_n,.)\big)(x_n)\\
&\le & P_{n,n+1}\big(\Phi(x_0,\ldots,x'_n,.)\big)(x'_n)
\end{eqnarray*}
where the first inequality follows from the non-negativity of the operator $P_{n,n+1}$ and the second follows from its monotony preserving property since $x_{n+1}\mapsto \Phi(x_{0:n+1})$ is non-decreasing. The function  $\Psi^{(n)}$, defined likewise, shares the same properties.
 
An induction assumption applied to the Markov chain $(X_k)_{0\le k\le n}$ completes the proof since
\begin{eqnarray*}
\E\Big(\Phi^{(n)}(X_{0:n})\Psi^{(n)}(X_{0:n})\Big)&\ge& \E \,\Phi^{(n)}(X_{0:n}) \E\, \Psi^{(n)}(X_{0:n})\\
&=&  \E \,\Phi(X_{0:n+1}) \E\, \Psi(X_{0:n+1}).\qquad \cqfd
\end{eqnarray*}

\ni $(b)$ By renumbering the $(n+1)$-tuple $(X_0,\ldots,X_n)$ we may assume $\tau  = id$. Then the transition $P_{k-1,k}(x_{k-1},dy)=\P_{X_k}(dy)$ does not depend  upon $x_{k-1}$ so that $P_{k-1,k}f$ is a constant function. Then~(\ref{MarkovfiniD}) holds as an equality and the monotony of $\Phi^{(n)}$ in each of its variable $x_0,\ldots,x_n$  is that of $\Phi$ for the same variables.  A careful inspection of the   proof  of claim~$(a)$ then shows that  the  weak co-monotony is enough to conclude.$\cqfd$

\bigskip 
\ni  {\sc Example.} Let $A, B\!\in {\cal B}or(\R^{n+1})$ be two Borel sets such that, for every $x=x_{0:n}\!\in A $, $x+t e_i\!\in A$ for every $t\!\in\R_+$ and every $i\!\in \{0,\ldots,n\}$ (where $e_i$ denotes the $i^{th}$ vector of the canonical basis of $\R^{n+1}$), {\em idem}  for $B$. Then for any $\R^d$-Markov chain $X=(X_k)_{0\le k\le n}$, having  monotony preserving transitions (in the sense of  Proposition~\ref{ProFinidim}$(a)$), we have
\[
\P\big((X_0,\ldots,X_n)\in A\cap B\big)\ge \P\big((X_0,\ldots,X_n)\in A\big)\P\big((X_0,\ldots,X_n)\in B\big).
\]

The monotony preserving property of the transitions $P_{k-1,k}$ cannot be relaxed as emphasized by the following easy counter-example.

\bigskip
\ni{\sc Counter-example.} Let  $X=(X_0,X_1)$ be a Gaussian bi-variate random vector with distribution $\displaystyle {\cal N}\Big(0; \Big[\begin{array}{cc}1&\rho\\ \rho&1\end{array}\Big]\Big)$ where the  correlation $\rho\!\in (-1,0)$). One checks that the transition $P_{0,1}(x_0,dx_1)$ reads on bounded or nonnegative Borel functions
\[
P_{0,1}(f)(x_0):= \E( f(X_1)\,|\,X_0=x_0) = \E f\Big(\rho\, x_0+\sqrt{1-\rho^2}\,Z\Big), \; Z\stackrel{\cal L}{\sim} {\cal N}(0;1).
\]
This shows that    $P_{0,1}$ is  monotony\dots {\em inverting}. In particular we have $\E\, X_0X_1 = \rho < 0 = \E\, X_0\,\E\, X_1$.
 In fact it is clear that $(X_0,X_1)$ satisfies the co-monotony principle if and only if $\rho\ge 0$.  In the next section we   extend this result to higher dimensional Gaussian vectors.

\subsection{More on the Gaussian case.} Let $X=(X_1,\ldots,X_d)$ be a centered Gaussian vector with covariance matrix $\Sigma= [\sigma_{ij}]_{1\le i,j\le d}$. 
This covariance matrix characterizes the distribution of $\Sigma$ so it characterizes as well  whether or not $X$ shares a co-monotony property in the sense of~(\ref{eq:comondimfin}). But can we {\em read easily } this  property on $\Sigma$? 

As mentioned above, a necessary condition for co-monotony is obviously that 
\[
\forall\, i,j\!\in\{1,\ldots,d\},\quad \sigma_{ij}= {\rm Cov}(X_i,X_j) \ge 0.
\]
In fact this simple condition does characterize co-monotony: this result, due to L. Pitt, is established in~\cite{PIT}.

\begin{Thm}[Pitt, 1982] A Gaussian random vector $X$ with covariance matrix $\Sigma= [\sigma_{ij}]_{1\le i,j\le d}$ satisfies a componentwise co-monotony principle
if and only if 
\[
\forall\, i,j\!\in\{1,\ldots,d\},\quad \sigma_{ij}  \ge 0.
\]
\end{Thm}

\noindent {\bf Remarks.} $\bullet$ Extensions have been proved in~\cite{JDPEPI}.  Typically, if $Z\sim{\cal N}(0;I_{d})$, under appropriate regularity and integrability assumptions on a function $h:\R^{d}\to \R$, one has 
$$ 
\Big(\forall\, x\!\in \R^{d},\quad \frac{\partial^2h}{\partial x_i\partial x_j} (x)\ge 0\Big)\Longrightarrow \Big(\sigma_{ij}\mapsto \E \big(h(\sqrt{\Sigma}Z)\big) \mbox{ is non-decreasing}\Big).
$$

\noindent $\bullet$ Another natural criterion for co-monotony --~theoretically straightforward although not easy to ``read" in practice  on the covariance matrix itself~-- is to make the assumption that there exists a matrix $A= [a_{ij}]_{1\le i\le d, 1\le j\le r}$, $r\!\in \N^*$,  
with {\em nonnegative entries} $a_{ij}\ge 0$ such that $\Sigma = AA^*$.  Then  $X\stackrel{\cal L}{\sim} AZ$, $Z\ \stackrel{\cal L}{\sim}{\cal N}(0;I_r)$. So every component is a linear combination with nonnegative coefficients of the components of $Z$ and Proposition~\ref{ProFinidim}$(b)$ straightforwardly implies that $X$ shares  the co-monotony property~(\ref{eq:comondimfin}).

However, surprisingly, this criterion is not a characterization in general: if $d\le 4$, symmetric matrices $\Sigma$ with nonnegative entries  can always be decomposed as $\Sigma= AA^*$ where $A$ has nonnegative entries. But if $d\ge  5$, this is no longer true. The negative answer is inspired by  a former counter-example -- originally due to Horn -- when $d=r \ge 5$, reported and justified in~\cite{HAL} (see Equations~(15.39) and~(15.53) and the lines that follow, see also~\cite{DIA} for an equivalent formulation). To be precise,  the nonnegative  symmetric $5\times 5$ matrix $\Sigma$ (with rank $4$) defined by
\[
\Sigma = \left[\begin{array}{ccccc}
1&0&0&1/2&1/2\\
0&1&3/4&0&1/2\\
0&3/4&1&1/2&0\\
1/2&0&1/2&1&0\\
1/2&1/2&0&0&1
\end{array}\right]
\]
 has    nonnegative entries but cannot be written $AA^*$ where $A$ has nonnegative entries. Another reference of interest about this question is~\cite{BAEM}, especially concerning the geometrical aspects of this problem.

\subsection{Application to the Euler scheme.} The Euler scheme of a  diffusion is an important example of  Markov chain to which one may wish to apply the co-monotony principle. 
Let $X=(X_t)_{t\in [0,T]}$ be a Brownian diffusion assumed to be solution to the stochastic differential equation
\[
SDE\equiv dX^x_t= b(t,X^x_t) dt+ \sigma(t,X^x_t) dW_t, \quad t\!\in [0,T],\quad X_0=x.
\]
Its Euler scheme with step $h=T/n$ and Brownian increments is entirely characterized by its transitions 
\[
P_{k,k+1}(f)(x) = \E\, f\Big(x+h b(t^n_k,x) +\sigma(t^n_k,x)\sqrt{h}\,Z\Big), \; Z\stackrel{{\cal L}}{\sim}{\cal N}(0;1), \;k=0,\ldots,n-1,
\]
where  $t^n_k=kh= \frac{k}{n}T$, $k=0,\ldots,n$.
One easily checks that if the function $b$ is Lipschitz continuous  in $x$ uniformly in $t\!\in [0,T]$ and if $\sigma(t,x)=\sigma(t)$ is deterministic and lies in $L^2([0,T],dt)$, then, for  large enough $n$, the  Euler transition $P$ is   monotony preserving. 

This follows from the fact that $x\mapsto x+hb(t,x)$ is non-decreasing provided  $h\!\in (0,\frac{1}{[b]_{\rm Lip}})$, where $[b]_{\rm Lip}$ is the uniform Lipschitz coefficient of $b$.
\section{Functional co-monotony principle}\label{Trois}
The aim of this section is to extend the above co-monotony principle to continuous time processes relying on the above multi-dimensional co-monotony result. To do so, we will view processes as random variables taking values in a path vector subspace $E\subset {\cal F}([0,T], \R)$ endowed with the (trace of the) Borel $\sigma$-field of pointwise convergence topology on $ {\cal F}([0,T], \R)$, namely $ \sigma\big(\pi_t, \, t\in [0,T]\big)$ where $\pi_t(\alpha)=\alpha(t)$, $\alpha\!\in E$. Consequently,  a process $X$ having $E$-valued paths can be seen as an $E$-valued random variable if and only if for every $t\!\in [0,T]$, $X_t$ is an $\R$-valued random variable (which is in some sense a tautology since it is  the lightest definition of a stochastic process).  

We consider on $E$ the   (partial) order  induced  by  the natural partial  ``pointwise order" on ${\cal F}([0,T], \R)$ defined by 
\[
\forall\, \alpha,\, \beta\!\in {\cal F}([0,T], \R),\; \; \alpha\le \beta\ \; \mbox{ if }\; \forall\, t\!\in [0,T],\;\alpha(t)\le \beta(t).
\]
\begin{Dfn} $(a)$ A measurable  functional $F:E\to \R$ is   {\em monotone} if it is either non-decreasing or non-increasing with respect to the order on $E$. 

\ss
\ni $(b)$ A {\em  pair}  of measurable  functionals are co-monotone if they are both monotone, with the same monotony. 
\end{Dfn}
Then, in order to establish a functional co-monotony principle (see the definition below)  our approach will be  transfer   a finite dimensional co-monotony principle satisfied by    appropriate converging (time) discretizations of the process $X$ of interest. Doing so we will need to  equip  $E$ with a topology ensuring the above convergence for the widest class of ($\P_{_X}$-$a.s.$ continuous) functionals. That is why we will   consider as often as  we can the sup-norm topology, not only    on   ${\cal C}([0,T], \R)$, but also on  the Skorokhod space $\D ([0,T], \R)$ of c\`adl\`ag (French acronym for right continuous left limited)  functions defined on $[0,T]$ since there are more continuous functionals for this topology than for   the  classical $J_1$ Skorokhod topology (having in mind that, furthermore,  $\D ([0,T], \R)$ is not a topological space for the lattest one). We recall that ${\cal D}_{_T}:= \sigma\big(\pi_t, \, t\in [0,T]\big)$ is the Borel $\sigma$-field related to both the $\|\,.\|_{\sup}$-norm and the $J_1$-topologies on the Skorokhod space. 

We will also consider (see Section~\ref{CadlagMarkov}) the space $L^p([0,T], \mu)$, $0<p<+\infty$, equipped  with its usual $L^p(\mu)$-norm where $\mu$ is a finite measure on $[0,T]$. In the latter case (which is not -- strictly speaking~-- a set of functions), we will consider the ``$\mu$-$a.e.$" (partial)  order 
\[
\alpha\le_{\mu} \beta\; \mbox{ if }\;\alpha(t)\le \beta(t) \; \mu(dt)\mbox{-}a.e.
\]
A functional $F$ which is monotone for the order $\le_{\mu}$ is called $\mu$-monotone. The definition of $\mu$-co-monotony follows likewise.

\medskip
A formal definition for the co-monotony  property on a partially ordered normed vector space  is the following.

 \begin{Dfn}\label{Def:2} A random variable $X$ whose paths take values in a partially ordered normed vector space $(E,\|\,.\, \|_E, \le_E)$~(\footnote{For every $\alpha$, $\beta$, $\gamma\!\in E$ and every $\lambda\!\in \R_+$, $\alpha\le \beta\Rightarrow \alpha+\gamma\le\beta+\gamma$ and $\lambda\ge 0\Rightarrow\lambda \alpha\le \lambda \beta$.}) satisfies a  co-monotony principle on $E$ if, for every bounded, co-monotone, $\P_{_X}$-$a.s.$ continuous, measurable functionals $F$, $G:E\to \R$,
\[
\E\, F(X)G(X)\ge \E\, F(X)\E\, G(X).
\]
\end{Dfn}

When $X$ is a stochastic process and $E$ is  its natural path space, we will often use the term {\em functional co-monotony principle}.

\bigskip
\noindent{\sc Extensions.} $\bullet$ The extension to square integrable or nonnegative $\P_{_X}$-$a.s.$ continuous functionals is canonical by a standard truncation procedure: replace $F$ by $F_N:=(-N)\vee \big(F\wedge N)$, $N>0$, and let $N$ go to infinity.

\smallskip\noindent  $\bullet$ More generally, the inequality also holds  for pairs of co-monotone functionals $F$, $G$ whose truncations $F_{_K}$ and $G_{_K}$ are  limits in $L^2(\P_{_X})$ of $\P_{_X}$-$a.s.$  continuous co-monotone functionals.

\subsection{A stability result for series of  independent random vectors}

We will  rely several times on the following proposition which shows that series of independent $E$-valued random vectors satisfying the co-monotony principle also share this property.

\begin{Pro}\label{LemTech} $(a)$ Let $(X_n)_{n\ge 1}$ be a sequence of independent $E$-valued random vectors defined on $(\Omega, {\cal A}, \P)$ where  $(E, \|\,.\,\|_E,\le)$ is a partially ordered normed vector space. Assume that, for every $n\ge 1$, $X_n$ satisfies a co-monotony principle on $E$.
Let  $(a_n)_{n\ge 1}$ be a sequence of real numbers such that the series $X= \sum_{n\ge 1} a_n X_n$ converges $a.s.$ for  the  norm $\|\,.\,\|_E$. Then $X$ satisfies a co-monotony principle on $E$.

\smallskip
\ni $(b)$ Assume furthermore that $(E, \|\,.\,\|_E)$ is a Banach space with an unconditional  norm, that the $X_n$ are nonnegative random vectors for the order on $E$ and that $\sum_{n\ge 1}X_n$ converges in $L_E^1(\P)$. Then, for every sequence of independent random variables $(A_n)_{n\ge 1}$ taking values in a fixed compact  interval of $\R$ and independent of $(X_n)_{n\ge 1}$,  the series $X= \sum_{n\ge 1}A_nX_n$ satisfies a co-monotony principle on $E$.
\end{Pro}

\ni{\bf Remark.} When $E= {\cal C}([0,T], \R)$,  L\'evy-It\^o-Nisio's Theorem (see $e.g.$~\cite{LETA}, Theorem 6.1, p.151) shows to some extent the equivalence between functional  convergence in distribution and $a.s.$ convergence for series of independent processes as above.

\bigskip
\ni {\bf Proof.} $(a)$ We may assume without loss of generality that the two functionals $F$ and $G$ are non-decreasing. We  first show the result for the sum of two independent processes, $i.e.$ we assume $a_k=0$, $k\ge 3$ (and $a_1a_2\neq 0$). 
%
%
By Fubini's Theorem
%
\[
\E\Big( F\big(a_1X_1+a_2X_2\big)G\big(a_1X_1+a_2X_2\big)\Big) = \E \Big( \left[\E \big( F\big(a_1 X_1+a_2 \alpha\big)G\big(a_1 X_1+a_2 \alpha\big)\big)\right ]_{|\alpha=X_2}\Big).
\]
Let Cont$(F)$ denote the set of elements of $E$ at which $F$ is continuous. It follows, still  from Fubini's Theorem, that 
\[
1=\P\big(a_1X_1+a_2X_2\!\in {\rm Cont}(F))\big)= \int\P_{X_2}(d\alpha_2)\P_{X_1}\big(F(a_1.+a_2\alpha_2)\big)
\]
so that $\P_{X_2}(d\alpha_2)$-$a.s.$ $\alpha_1\mapsto F(a_1\alpha_1+a_2\alpha_2\big)$ and $\alpha_1\mapsto G(a_1\alpha_1+a_2\alpha_2\big)$ are $\P_{X_1}$-$a.s.$ continuous. Noting that these functionals are co-monotone (non-decreasing if $a_1\ge 0$, non-increasing if $a_1\le 0$), this implies
\[
\E\big( F(a_1X_1+a_2\alpha_2)G(a_1X_1+a_2\alpha_2)\big)\ge \E\big( F(a_1X_1+a_2\alpha_2)\big)\E\big(G(a_1X_1+a_2\alpha_2)\big).
\]


Now, both
$$
\alpha_2\mapsto \E\,F(a_1X_1+a_2\alpha_2)=\int\P_{X_1}(d\alpha_1) F(a_1\alpha_1+a_2\alpha_2)
$$ 
and  
$$
\alpha_2\mapsto \E G(a_1X_1+a_2\alpha_2)=\int\P_{X_1}(d\alpha_1) G(a_1\alpha_1+a_2\alpha_2)
$$ 
are co-monotone (non-decreasing if $a_2\ge 0$, non-increasing if $a_2\le 0$) and one checks that both are $\P_{X_2}(d\alpha_2)$-$a.s.$ continuous which implies in turn that
\begin{eqnarray*}
\E\Big(\!\big[ \E\, F(a_1X_1+a_2\,\alpha)\big]_{|\alpha=X_2} \!\!\!&\!\!&\!\! \!\!\!\!\!\!\big[\E\,G(a_1X_1+a_2\,\alpha)\big]_{|\alpha=X_2}\!\Big)\\
\!\!\!&\!\ge\!&
\E\Big(\!\big[ \E\, F(a_1X_1+a_2\,\alpha)\big]_{|\alpha=X_2}
\Big)\E\Big(\big[\E\,G(a_1X_1+a_2\alpha)\big]_{|\alpha=X_2}\Big)\\
\!\!\!&\!=\!&\E\Big(F(a_1X_1+a_2X_2)\Big)\E\Big(G(a_1X_1+a_2X_2)\Big)
\end{eqnarray*}
where we used again Fubini's Theorem in the second line.

One extends this result by induction to the case where $X=X_1+\dots+X_n$.

\smallskip To make $n$ go to infinity, we proceed as follows: let ${\cal G}_n=\sigma(X_k,\, k\ge n+1)$. By the reverse martingale convergence theorem, we know that 
for any bounded measurable functional $\Phi:E\to \R$,
\[
\E\,( \Phi(X)\,|\, {\cal G}_n) = \left[\E\,\Phi(X_1+\cdots+X_n +\widetilde \alpha_n)\right]_{\widetilde \alpha_n =\widetilde X_n}
\]
where $\widetilde X_n =X-(X_1+\cdots+X_n)$. We know from the above case $n=2$ that, for $\Phi=F$ or  $G$,   one has $\P_{\widetilde X_n}(d\widetilde \alpha_n)$-$a.s.$, $\alpha\mapsto \Phi(\alpha+\widetilde \alpha_n)$ is $\P_{X_1+\cdots+X_n}(d\alpha)$-continuous so that
\[
\E\,FG(X_1+\cdots+X_n +\widetilde \alpha_n) \ge \E\,F(X_1+\cdots+X_n +\widetilde \alpha_n)\E\,G(X_1+\cdots+X_n +\widetilde \alpha_n).
\]
This equality also reads
\[
\E\,( FG(X)\,|\, {\cal G}_n) \ge \E\,( F(X)\,|\, {\cal G}_n) \E\,( G(X)\,|\, {\cal G}_n) 
\]
which in turn implies by letting $n\to \infty$ 
\[
\E\,F(X)G(X)\ge \E\, F(X)\E\,G(X)
\]
owing to the reverse martingale convergence theorem.

\smallskip
\ni $(b)$ For every bounded sequence $(a_n)_{n\ge 1}$, it follows from  the unconditionality of the norm $\|\,.\,\|_E$ that $\sum_{n\ge 1}a_nX_n$ $a.s.$ converges in $L^1_E(\P)$. Then it follows from   $(a)$ that, for every $n\ge 1$,
\[
\E\left(F\Big(\sum_{k= 1}^na_k X_k\Big)G\Big(\sum_{k= 1}^na_k X_k\Big)\right) \ge \E\, F\Big(\sum_{k= 1}^na_k X_k\Big)\E\, G\Big(\sum_{k= 1}^na_k X_k\Big).
\]
Now for every $k\!\in \{1,\ldots,n\}$, the function defined on the real line by $a_k\mapsto    \E\, F\big(\sum_{i= 1}^na_i X_i\big)$ has the same monotony as $F$ since $X_k\ge 0$ and is bounded. Consequently for any pair $F$, $G$ of bounded co-monotone Borel functionals,
\begin{eqnarray*}
\E\Big(\big[\E\, F(\sum_{k= 1}^na_k X_k)\big]_{|a_{1: n}=A_{1: n}} \!\!\!&\!\!&\!\! \!\!\!\!\!\!\!\!\!\big[\E\, G(\sum_{k= 1}^na_k X_k)\big]_{|a_{1: n}=A_{1: n}}\Big)\\&\ge& \E\Big(\big[\E\, F(\sum_{k= 1}^na_k X_k)\big]_{|a_{1: n}=A_{1: n}}  \Big)\times    \E\Big(\big[\E\, G(\sum_{k= 1}^na_k X_k)\big]_{|a_{1: n}=A_{1:n}}\Big).
\end{eqnarray*}
The conclusion follows  for a fixed $n\ge 1$ by preconditioning. One concludes by letting $n$ go to infinity since $F$ and $G$ are continuous.~$\cqfd$

\bigskip
 \noindent  
 {\bf A first application to Gaussian processes.} Let $X=(X_t)_{t\in[0,T]}$ be a continuous centered Gaussian process with a covariance operator $C_X$ defined on the Hilbert space $L^2_{_T}:=L^2([0,T],dt)$ into itself
 by 
 \[
 C_X(f) = \E( \langle f,X\rangle_{L^2_T}X)= \int_0^T \E(X_sX_.)f(s)ds\!\in L^2_{_T}.
 \]
The process $X$ can be seen as a random vector taking values in the separable Banach space ${\cal C}([0,T], \R)$ (equipped with the $\sup$-norm). Assume that $C_X$ admits a   decomposition as follows 
 \[
 C_X = AA^*, \qquad A: (K, |\,.\,|_K)\longrightarrow {\cal C}([0,T], \R), \, A\; \mbox{continuous linear mapping},
 \]
where $(K, |\,.\,|_K)$ is a separable Hilbert space.

Then, we know from Proposition~1 (and Theorem~1) in~\cite{LUPA5}, that for any orthonormal basis (or even any {\em Parseval frame}, see~\cite{LUPA5}) $(e_n)_{n\ge 1}$ of $K$ that the sequence $(A(e_n))_{n\ge 1}$ is {\em admissible} for the process $X$ in the following sense: for any i.i.d. sequence $(\xi_n)_{n\ge 1}$ of normally distributed random variables defined on a probability space $(\Omega, {\cal A}, \P)$
 \[
 \left\{\begin{array}{ll}  (i)&  \displaystyle \sum_{n\ge 1} \xi_n A(e_n)\mbox{ $a.s.$ converges   in $({\cal C}([0,T], \R), \|\,.\,\|_{\sup})$}\\
 (ii) &  \displaystyle\sum_{n\ge 1} \xi_n A(e_n) \stackrel{\cal L}{\sim} X.\end{array}\right.
 \]
 Assume furthermore that {\em  all the continuous functions $A(e_n)$ are  nonnegative}. Then,  for every $n\ge 1$, the continuous stochastic process $X_n = \xi_n A(e_n)$ satisfies a co-monotony principle (for the natural pointwise partial order on ${\cal C}([0,T],\R)$). This makes up a sequence of independent random elements of ${\cal C}([0,T], \R)$. It follows from Proposition~\ref{LemTech}$(a)$ that the process $X$ satisfies a co-monotony principle.

  \bigskip
 \noindent {\sc Example:} Let us  consider the standard Brownian motion $W$ with covariance function $\E \,W_tW_s= s\wedge t$. One checks that $C_W = AA^*$ where $A: L^2_{_T}\to {\cal C}([0,T],\R)$ is defined by 
 \[
 A f \equiv\Big(t\mapsto \int_0^t f(s)ds\Big)\!\in  {\cal C}([0,T],\R).
 \]
Applied to the orthonormal basis $e_n(t) = \sqrt{\frac 2T} \sin \big(\pi n \frac tT\big)$, $n\ge1$, we get 
$$
\displaystyle A(e_n) (t)=  \sqrt{2T}\,\frac{1-\cos\big(\frac{\pi n t}{T} \big)}{\pi n}\ge 0, \;t\in [0,T],\;n\ge 1,
$$ 
so that  
\[
  \widetilde W=  \sqrt{2T} \sum_{n\ge 1}\frac{\xi_n}{\pi n}\Big(1-\cos\big(\pi n \frac .T\big)\Big),\; (\xi_n)_{n\ge 1} \mbox{ i.i.d.}, \; \xi_1\stackrel{\cal L}{\sim} {\cal N}(0;1),
\] 
is an $a.s.$ converging series for the $\sup$-norm which defines  a standard Brownian motion. As a consequence, the standard Brownian motion  satisfies a co-monotony principle (in the sense of Definition~\ref{Def:2}).~(\footnote{The fact that $A(L^2_{_T})$ is the Cameron-Martin space $i.e.$ the reproducing space of the covariance operator, which is obvious here, is a general fact for any such decomposition (see~\cite{LUPA5}).})
 
 \medskip
 As we will see further on in Section~\ref{GaussProc}, the above approach is clearly neither the most elementary way nor the most straightforward  to establish the co-monotony principle for the Wiener process.  Furthermore, the above criterion is not an equivalence as emphasized in a finite dimensional setting: a continuous Gaussian process $X$ may satisfy a functional co-monotony principle albeit its covariance  operator $C_X$ admits representation of the form  $C_X=AA^*$ for which there exists  an orthonormal basis (or even Parseval frame, see~\cite{LUPA5}) whose image by $A$ is made of nonnegative functions. Thus, no such decomposition is known to us for the fractional Brownian motion (with Hurst constant $H\neq \frac 12$) although it satisfies a co-monotony principle (see Section~\ref{fracW} further on).

\subsection{From $[0,T]$ to $\R_+$}
We state our results  for processes defined on a finite interval $[0,T]$. However they can be extended canonically on $\R_+$, provided that there exists a sequence of positive real constants $T_N\uparrow+ \infty$ such that $\P_X(d\alpha)$-$a.s.$ on $E\subset {\cal F}(\R_+, \R)$,  $\alpha^{T_n}\equiv \big(t\mapsto \alpha(t\wedge T_n)\big)$ converges in $E$ toward $\alpha$ for the topology on $E$. Such a sequence does exist the topology of convergence on compact sets but also for the Skorokhod topology on the positive real line. Then the transfer of co- and anti-monotony property  (if any)  from the stopped process $X^{T_n}$ to $X$ is straightforward for bounded functionals. The extension to square integrable or nonnegative functionals follows by the usual truncation arguments.

\section{Application to pathwise continuous processes}

 In this section functional co-monotony principle is always considered {\em on the normed vector space $({\cal C}([0,T],\R), \|\,.\,\|_{\sup})$}. We  will use implicitly that its  Borel $\sigma$-field of   is $\sigma\big(\pi_t, \!\in [0,T]\big)$, where $\pi_t(\alpha)=\alpha(t)$ for every $\alpha\!\in {\cal C}([0,T], \R)$ and every $t\!\in[0,T]$  (see~\cite{BIL}, Chapter~2).

\begin{Pro}\label{ContProc}Let $X= (X_t)_{t\in [0,T]}$ be a pathwise continuous  process defined on a probability space $(\Omega,{\cal A},\P)$ sharing the following finite dimensional co-monotony property: for every integer $n\ge 1$ and every subdivision $(t_1,\ldots, t_n)\!\in [0,T]^n$, $0\le t_1<t_2<\cdots<t_n \le T$, the random vector $(X_{t_k})_{1\le k\le n}$ satisfies a componentwise co-monotony principle. Then $X$ satisfies a functional co-monotony principle on its path space ${\cal C}([0,T], \R)$.
\end{Pro}

\ni {\bf Proof.} We may assume that $F$ and $G$ are both non-decreasing for the natural order on $ {\cal C}([0,T], \R)$. Let $n\in \N$, $n\ge 1$. We introduce the uniform subdivision $t^n_k =\frac{kT}{n}$, $k=0,\ldots,n$ and, for every function $\alpha\!\in {\cal C}([0,T], \R)$,  the canonical linear interpolation approximation
\[
\alpha^{(n)}(t) = \frac{t^n_{k+1}-t}{t^n_{k+1}-t^n_k}\alpha(t^n_k) + \frac{t-t^n_k}{t^n_{k+1}-t^n_k}\alpha(t^n_{k+1}),\; t\!\in [t^n_k,t^n_{k+1}],\; k=0,\ldots,n-1.
\]
  One checks that $\|\alpha-\alpha^{(n)}\|_{\sup}\le w(\alpha, T/n)$ goes to $0$ as $n\to \infty$ where $w(\alpha, .)$ denotes the uniform continuity modulus of $\alpha$. As a consequence, $X$ having $a.s.$ continuous paths by assumption, the sequence $(X^{(n)})_{\ge 1} $ of interpolations of $X$ $a.s.$ uniformly converges toward $X$.

  Then set for every $n\ge 1$ and every $x=x_{0:n}\!\in \R^{n+1}$,
  \[
  \chi^n(x,t) = \frac{t^n_{k+1}-t}{t^n_{k+1}-t^n_k}x_k + \frac{t-t^n_k}{t^n_{k+1}-t^n_k}x_{k+1},\; t\!\in [t^n_k,t^n_{k+1}],\; k=0,\ldots,n-1
  \]
 and  
 \[
 F_n(x) =F(\chi^n(x,.)).
  \]
  It is clear that if $x\le x'$ in $\R^{n+1}$ (in the componentwise sense) then $F_n(x)\le F_n(x')$ since $\chi^n(x,.)\le \chi^n(x',.)$ as functions. This is equivalent to the fact that $F_n$ is non-decreasing in each of its variables. 
  
  On the other hand $X^{(n)}= \chi\big((X_{t^n_k})_{0\le k\le n},.\big)$ so that $F_n\big((X_{t^n_k})_{0\le k\le n}\big)=F(X^{(n)})$. The sequence $(X_{t^n_k})_{0\le k\le n}$ satisfies a componentwise co-monotony principle. As a consequence, it follows from Proposition~\ref{ProFinidim} that if $F$ and $G$ are bounded, for every $n\!\in \N$,
  \[
  \E \,F(X^{(n)})G(X^{(n)})\ge \E\,F(X^{(n)}) \E  \,G(X^{(n)}).
  \]
One derives the expected inequality by letting $n$ go to infinity since $F$ and $G$ are continuous with respect to the $\sup$-norm. The extension to unbounded functionals $F$ or $G$ follows by the usual truncation arguments.~$\cqfd$
  
%

\subsection{Continuous Gaussian processes}\label{GaussProc}
Let $X=(X_t)_{t\in [0,T]}$ be a continuous centered Gaussian process. Its (continuous) covariance function $C_{_X}$ is defined on $[0,T]^2$ as follows 
\[
\forall\, s,\, t\!\in [0,T],\qquad C_{_X}(s,t) = \E\, X_sX_t.
\]
We establish below the functional  counterpart of Pitt's Theorem.

\begin{Thm}[Functional Pitt's Theorem] A continuous  Gaussian process  $X=(X_t)_{t\in [0,T]}$  with  a covariance operator $C_{_X}$ satisfies a functional co-monotony principle if and only if, for every $s,\, t\!\in [0,T],\;C_{_X}(s,t)\ge 0$.
\end{Thm} 

\noindent{\bf Proof.} 
For every $n\ge 1$, $(X^n_{t^n_k})_{0\le k\le n}$ is a  Gaussian random vector satisfying a componentwise co-monotony principle since its covariance matrix $\Sigma^n = \big[C_{_X}(t^n_k,t^n_{\ell})\big]_{0\le k,\ell\le n}$ having nonnegative entries. 
One concludes by the above Proposition~\ref{ContProc}.$\cqfd$

\bigskip
\noindent We briefly  inspect below several classical classes of Gaussian processes.

\subsubsection{Brownian motion, Brownian bridge, Wiener integrals} The covariance of the standard Brownian motion $W$ is given by $C_W(s,t) = s\wedge t \ge 0$, $s$, $t\!\in [0,T]$ and that of the Brownian bridge over $[0,T]$ defined by $X_t = W_t-\frac tTW_T$, $t\!\in [0,T]$, is given by  for every $s$, $t\!\in [0,T]$  by 
$$
C_X(s,t) = s\wedge t-\frac{st}{T}\ge 0.
$$

As for Wiener integrals, let $X_t= \displaystyle \int_0^t f(s)dW_s$, $t\!\in [0,T]$, where $f\!\in L^2([0,T],dt)$. The process  $X$ admits a continuous modification and its covariance function is given by  $C_X(s,t) = \displaystyle \int_{s\wedge t}^{s\vee t} f^2(u)du\ge 0$.


\subsubsection{Liouville processes}

\begin{Dfn} Let $f:[0,T]\to \R$ be a locally $\rho$-H\"older function, $\rho\!\in (0,1]$, in the following sense: there exists $\varphi\!\in {\cal L}^2([0,T],dt)$, $\rho\!\in (0,1]$, $a\!\in(0,+\infty)$  such that 
\begin{equation}
(L_{\rho,a})\equiv \left\{\begin{array}{ll}
(i)&  \forall\, t,\,t'\!\in [0,T],\quad |f(t)-f(t')|\le [f]_{\rho,\varphi}|t-t|^{\rho}\varphi(t\wedge t')\\
&\\
(ii)&  \displaystyle \int_0^tf^2(s)ds=O(t^a).
\end{array}\right.\hskip 2 cm 
\end{equation}
Then the Gaussian  process defined for every $t\!\in[0,T]$ by 
\[
X_t =\int_0^t f(t-s)dW_s
\]
admits a continuous modification  called {\em Liouville process} (related to $f$) with covariance function
\[
C_{_X}(s,t)= \int_0^{s\wedge t}f(t-u)f(s-u)du.
\] 
\end{Dfn}
  
\ni  {\sc Justification:} First note that $f\!\in {\cal L}^2(dt)$ since $f(t)|\le |f(0)| +t^{\rho}|\varphi(0)|$. Then, for every $t,t'\!\in [0,T]$,  $t\le t'$, 
  \[
  X_{t'}-X_t =\int_{t}^{t'} f(t'-s)dW_s +\int_0^t(f(t'-s)-f(t-s))dW_s
  \]
  so that 
  \begin{eqnarray*}
  \E| X_{t'}-X_t|^2 &= &\int_t^{t'}f^2(t'-s)ds +\int_0^t \big(f(t'-s)-f(t-s)\big)^2ds\\
  &\le& \int_0^{t'-t}f^2(s)ds+  [f]^2_{\rho,\varphi} |t'-t|^{2\rho}\int_0^T\varphi^2(s)ds\\
  &\le& C_f(|t'-t|^a+ |t'-t|^{2\rho} ) \\
  &\le & C_f|t'-t|^{(2\rho)\wedge a}
  \end{eqnarray*}
  so that, using the Gaussian feature of the process $X$
 \[
  \E| X_{t'}-X_t|^p\le C_{f,p}|t'-t|^{p(\rho\wedge a/2)}
 \]
  for every $p\ge 2$  which in turn implies, owing to Kolmogorov's continuity criterion,  that $(X_t)_{t\in[0,T]}$ admits a version  $\big(\rho\wedge \frac a2)-\eta\big)$-H\" older continuous for any small enough $\eta >0$.
  
  \begin{Pro}\label{Liouville} Let $\displaystyle X_t=\int_0^t \!f(t-s)dW_s$, $t\!\in [0,T]$, be  a continuous Liouville process where $W$ is a standard B.M. defined on a probability space $(\Omega,{\cal A}, \P)$ and  $f$  satisfies $(L_{\rho,a})$  for a couple $(\rho,a)\!\in (0,1]\times (0,+\infty)$. If furthermore $f$ is $\lambda$-$a.e.$ nonnegative, then $X$ satisfies a co-monotony principle.  
 \end{Pro}

 The proof of the proposition is straightforward since $f\ge 0$ $\lambda$-$a.e.$ implies that the covariance function $C_{_X}$ is nonnegative.
 
 \bigskip
  \ni {{\sc Example:} If $f(u)= u^{H-1/2}$ with $H\!\!\in (0,1]$, then $f$ satisfies $(L_{a,\rho})$ with  $a=2H$ and $\rho \!\in (0, \frac 12 -H)$ if $H< \frac12$, $\rho =H-\frac 12$ if $H > \frac12$ (and $\rho=1$ if $H=\frac 12$). This corresponds to the pseudo-fractional Brownian motion with Hurst constant $H$.

\subsubsection{Wiener integrals depending on a parameter} Now we consider a  class of processes which is wider  than  Liouville's class and for which we provide a slightly less refined criterion of existence (as a pathwise continuous process).
\[
X_t = \int_{0}^{\infty}f(t,s)dW_s,\quad t\!\in [0,T], 
\]
where $f:[0,T]\times \R_+ \to \R$ satisfies a dominated $\rho$-H\"older assumption reading as  follows: there exists $\varphi\!\in{\cal  L}^2(\R,dt)$, non-increasing, and $\rho\!\in (0,1]$ such that 
\[
(L'_{\rho})\;\equiv \;  \forall\, t,\,t'\!\in [0,T], \; \forall\, s\!\in \R_+,\quad |f(t',s)-f(t,s)|\le [f]_{\rho} |t'-t|^{\rho}\varphi(s).\hskip 3cm
\]
Such a process has a continuous modification since $t\mapsto X_t$  is $\rho$-H\" older from $[0,T]$ to $L^2(\P)$ and Gaussian (still owing to Kolmogorov's continuity criterion). 



\smallskip
As for Liouville processes, if furthermore, for every $t\!\in [0,T]$,  {\em $f(t,.)$ is $\lambda_1$-$a.e.$ nonnegative, then the process $X$ satisfies a co-monotony principle}. 

\bigskip
\ni {\sc Example.} Let $f_H(t,s) =(t+s)^{H-\frac 12}-s^{H-\frac 12}$, $H\!\in (0,1]$. 

\begin{itemize} \item If $H\ge 1/2$, $f_H$ satisfies $(L'_{H-\frac 12})$  with $\varphi(s) = 1$. 
\item If $H\!\in (0,\frac 12]$, $f_H$ satisfies $(L'_{\frac 12-H})$ with $\varphi(s) = s^{2H-1}$.
\end{itemize}
 \subsubsection{Fractional Brownian motion with Hurst constant $H\!\in (0,1]$}\label{fracW}
 The fractional Brownian motion is a continuous Gaussian process  characterized by its covariance function defined by
 \[
 \forall\, s,\, t\!\in [0,T], \quad C^H(s,t)= \frac 12 \big(t^{2H}+s^{2H}-|t-s|^{2H}\big).
 \]
 
 Since $u\mapsto u^H$ is $H$-H\"older and $|t-s|^2\le t^2+s^2$, it is clear that $C(s,t)\ge 0$. Consequently, {\em  the fractional Brownian motion satisfies a co-monotony principle}.

\bigskip
\noindent {\bf Remark.}  An alternative approach could be to rely on  the celebrated Mandelbrot-Van Ness representation of the fractional Brownian motion with Hurst constant $H\!\in (0,1]$, given by
 \[
 B^H_t= \underbrace{\int_0^{+\infty} \big((t+s)^{H-\frac 12}-s^{H-\frac 12} \big)dW^1_s}_{_{=:B^{H,1}_t}}+ \underbrace{\int_0^t |t-s|^{H-\frac 12} dW^2_s}_{=:B^{H,2}_t} 
 \]
 where $W^1$ and $W^2$ are independent standard Brownian motions. These two Wiener integrals define pathwise continuous independent processes, both satisfying the co-monotony principle for $\P_{B^{H,i}}$-$a.s.$ $\|\,.\,\|_{\sup}$-continuous functionals, consequently their sum satisfies a co-monotony principle for $\P_{B^{H}}$-$a.s.$ $\|\,.\,\|_{\sup}$-continuous functionals owing to Lemma~\ref{LemTech}$(b)$.

\subsection{Continuous Markov Processes, Brownian diffusion processes} \label{ContMarkov}

\begin{Pro}\label{ContMarkov} Let $X= (X_t)_{t\in [0,T]}$ be a pathwise continuous Markov process defined on a probability space $(\Omega,{\cal A},\P)$ with  transition operators $(P_{s,t})_{0\le s\le t\le T}$ satisfying the monotony preserving property. Then $X$ satisfies a functional co-monotony principle on its path space ${\cal C}([0,T], \R)$.
\end{Pro}

\ni {\bf Proof.} The sequence $(X_{t^n_k})_{0\le k\le n}$ is a discrete time Markov chain whose transition operators $P_{t^n_k,t^n_{k+1}}$, $k=0,\ldots,n-1$  satisfy the monotony preserving property. ~$\cqfd$

\bigskip
  The main application of this result  is the {\em co-monotony principle for Brownian diffusions} ($i.e.$  solutions of stochastic differential equations driven by a standard Brownian motion). We consider the Brownian diffusion
\[
(SDE)\; \equiv \; dX^x_t= b(t,X^x_t)dt +\sigma(t,X^x_t) dW_t, \; X^x_0=x  
\]
where $b: [0,T]\times \R\to \R$ is continuous, Lipschitz continuous  in $x$, uniformly in $t\!\in [0,T]$ and  $\sigma: [0,T]\times \R\to \R$ (is continuous) with linear growth in $x$, uniformly in $t\!\in [0,T]$ and satisfies
\[
\forall\, x,\, y\!\in \R ,\; \forall\,t\!\in [0,T], \quad |\sigma(t,x)-\sigma(t,y)|\le \rho(|x-y|),\; 
\]
where 
\[
\rho:\R\to\R \mbox{ is increasing}, \;\rho(0)=0 \mbox{ and } \int_{0+}\frac{du}{\rho^2(u)}=+\infty.
\] 
Then  the equation $(SDE)$ satisfies a weak existence property since $b$ and $\sigma$ are continuous with linear growth: its Euler scheme weakly functionally converges  to a weak solution of $(SDE)$ as its step  $T/n$  goes to $0$ (see Theorem~5.3 in~\cite{JASH}). It also satisfies  a strong  uniqueness property (see Proposition~2.13 in~\cite{KASH}) hence a strong existence-uniqueness property. This implies the existence of (Feller) Markov transitions $(P_{s,t}(x,dy))_{t\ge s\ge 0}$  such that, $a.s.$ for every $x\!\in\R$,  $P_{s,t}(f)(x) = \E( f(X_t)\,|\,X_s=x)$ (see $e.g.$ Theorem~1.9 in~\cite{REYO}).  Furthermore the flow $(X^x_t)_{x\in \R, t\in [0,T]}$ satisfies a comparison principle (Yamada-Watanabe's Theorem)  $i.e.$   for every $x, \,x'\!\in \R$, $x\le x'$, $\P$-$a.s.$, for every $t\!\in [0,T]$, $X_t^x\le X^{x'}_t$.  The functional co-monotony principle follows immediately  since it implies that the Markov transitions $P_{s,t}$ are monotony preserving.


\bs
\ni {\bf Remarks.} $\bullet$ It is to be noticed that the Euler scheme does not share the componentwise co-monotony principle as a Markov chain in full generality, in particular, when $\sigma$ does depend on $x$.  So, the above result for diffusions is not a simple transfer from the discrete time case. However when $\sigma(x,t)=\sigma(t)$ with  $\sigma:[0,T] \to \R$  continuous,  the result can be transferred from   the Euler scheme (see~Section~\ref{Deux}) since this scheme functionally weakly  converges toward  $X^x$ as its step $T/n$ goes to $0$.

\medskip
\ni $\bullet$ This result is strongly related to strong uniqueness of solutions of $(SDE)$. The above conditions are not minimal, see $e.g.$~\cite{ENG} for more insights on these aspects.

%


\section{Functional co-monotony principle for c\`adl\`ag processes}\label{Quatre}

The most natural idea is to mimick Proposition~\ref{ContProc} by simply replacing $({\cal C}([0,T],\R), \|\,.\,\|_{\sup})$ by the space $I\!\!D([0,T],\R)$ of c\`adl\`ag functions endowed with the $J_1$  Skorokhod toplogy (see~\cite{BIL}, Chapter 3). Although this is not a topological vector space (which make some results fail like Proposition~\ref{ContProc}), this approach yields non-trivial results. To be precise, let us consider, instead of the interpolation operator of a continuous  function on the uniform  subdivision $(t^n_k)_{0\le k\le n}$, the {\em stepwise constant approximation} operator define on every function $\alpha\!\in I\!\!D([0,T],\R)$ by

\begin{equation}\label{stepwise}
\widetilde \alpha^{(n)}=  \sum_{k=1}^n \alpha(t^n_{k-1})\mbox{\bf 1}_{[t^n_{k-1},t^n_k)} +\alpha(T)\mbox{\bf 1}_{\{T\}}, \quad t^n_k=\frac {kT}{n}, \; k=0,\ldots,n.
\end{equation}

It follows from Proposition 6.37 in~\cite{JASH}, Chapter VI (second edition), that $\alpha^{(n)}\to \alpha$ for the Skorokhod topology. Then by simply mimicking the proof of Proposition~\ref{ContProc}, we get the following result.

\begin{Pro}Let $X= (X_t)_{t\in [0,T]}$ be a pathwise c\`adl\`ag  process defined on a probability space $(\Omega,{\cal A},\P)$ sharing the  finite dimensional co-monotony property. Then $X$ satisfies a functional co-monotony principle on its path space $I\!\!D([0,T], \R)$ for the $J_1$ Skorokhod topology.
\end{Pro}

Thus if one considers now a general c\`adl\`ag process with independent increments (PII) $(X_t)_{t\ge 0}$ defined on a probability space $(\Omega,{\cal A}, \P)$, it is clear, $e.g.$  from Proposition~\ref{ProFinidim}$(a)$, that $X$ shares the finite dimensional co-monotony property since $(X_{t^n_k})_{0\le k\le n}$ is a Markov chain whose transitions $P_{t^n_{k-1},t^n_k}(x_{k-1},dy) = {\cal L}\big(x_{k-1} + X_{t^n_k}-X_{t^n_{k-1}}\big)$, $k=1,\ldots,n$, are clearly monotony preserving. 

\medskip  Consequently {\em any c\`adl\`ag process with independent increments (PII) $(X_t)_{t\ge 0}$ satisfies a functional co-monotony principle on  its path space $I\!\!D([0,T], \R)$ for the $J_1$ Skorokhod topology}.

\medskip Note that a continuous process which satisfies the functional co-monotony property on its path space for the sup-norm will always satisfy the functional co-monotony principle on $I\!\!D([0,T], \R)$ for the $J_1$ topology since, when $\alpha$ is continuous, convergence of a sequence $(\alpha_n)$ to $\alpha$  for the  sup-norm   and the Skorokhod topology coincide.

\medskip However this result is not fully satisfactory since there are not so many functionals which are continuous or even $\P_{_X}$-$a.s.$ continuous with respect to the Skorokhod topology. Thus the partial maxima functional $\alpha \mapsto \sup_{s\in [0,t]}|\alpha(s)|$ is not Skorokhod continuous if $\alpha $ is not continuous at $t$  (except when $t=T$) and this functional co-monotony principle will fail, $e.g.$  for any  process $X$ having a fixed discontinuity at $t$. This is the reason why we establish in the next subsection a functional co-monotony principle for general {\em PII} on $I\!\!D([0,T], \R)$ {\em endowed with the sup-norm topology}.

\subsection{Sup-norm co-monotony for processes with independent increments}
We consider a general c\`adl\`ag process with independent increments ({\em PII}) $(X_t)_{t\ge 0}$ defined on a probability space $(\Omega,{\cal A}, \P)$. We rely on its L\'evy-Khintchine decomposition  as exposed in~\cite{JASH}, chap. II, section~3. First one can decompose $X$ as the sum 
\[
X= X^{(1)}\stackrel{\perp\!\!\!\perp}{+}X^{(2)}
\]
where $X^{(1)}$ and $X^{(2)}$ are two independent PII:  $X^{(1)}$  is a {\em PII} {\em without fixed discontinuities} and  $X^{(2)}$ is a   pure jump PII,  possibly jumping only at  a deterministic sequence of  times, namely 
\[
X^{(2)}_t =\sum_{n\ge 1}U_n\mbox{\bf 1}_{\{t_n\le t\}}, \;t\!\in [0,T],
\]
where $(t_n)_{n\ge 1}$ is a sequence of $[0,T]$-valued real numbers and $(U_n)_{n\ge }$ is a sequence of independent random variables satisfying the usual assumption of the three series theorem
\begin{eqnarray*}
&&\sum_n \P(|U_n|\ge 1)<+\infty,\\
&& \sum_{n} \E U_n \mbox{\bf 1}_ {\{|U_n|\le 1\}}<+\infty,\\ 
&& \sum_{n} \E \big(U^2_n \mbox{\bf 1}_ {\{|U_n|\le 1\}}-(\E\, U_n \mbox{\bf 1}_ {\{|U_n|\le 1\}})^2\big)<+\infty.
\end{eqnarray*}

\begin{Pro} A c\`adl\`ag {\em PII} satisfies a co-monotony  principle on $(\D([0,T], \R), \|\,.\,\|_{\sup})$.
\end{Pro}

\ni {\bf Proof.} Owing to Lemma~\ref{LemTech}$(a)$, we will inspect successively the cases of {\em PII} without fixed discontinuities and of pure jumps.

\medskip
\ni{\sc Step~1.} {\em $X$ is a {\em PII} without fixed discontinuities}: This means that $X^{(2)}\equiv 0$.  The classical (pathwise) L\'evy-Khintchine formula for {\em PII} without fixed discontinuities says that, a truncation level $\varepsilon>0$ being fixed, $X$ reads as the  sum of three mutually independent processes as follows
\[
\forall\, t\!\in [0,T], \qquad X_t = b^{\varepsilon}(t) + W_{c(t)} \stackrel{\perp\!\!\!\perp}{+} \sum_{s\le t } \Delta X_s \mbox{\bf 1}_{\{|\Delta X_s|> \varepsilon\}} \stackrel{\perp\!\!\!\perp}{+} M^{\varepsilon}_t
\]
where $b^{\varepsilon}$ is a  continuous function on $[0,T]$, $c$ is a nonnegative non-decreasing continuous function on $[0,T]$ with $c(0)=0$ and $M^{\varepsilon}_t$ is a pure jump martingale satisfying
\[
\E \left(\sup_{s\in [0,t]}|M^{\varepsilon}_s|^2\right)\le 4 \int_{\R\setminus\{0\}} x^2\mbox{\bf 1}_{\{|x|\le \varepsilon\}} \nu^X ([0,t]\times dx)
\]
 where the measure $\nu^X$ is the L\'evy measure of $X$, $i.e.$ the dual predictable projection of the jump measure $\mu^X(ds,dx)=\sum_{s\in [0,T]} \mbox{\bf 1}_{\{\Delta X_s\neq 0\}}\Delta X_s$. The L\'evy measure is characterized by the fact  that, for every bounded Borel function $g:\R\to\R$  null in the neighbourhood of $0$, 
\[
\left(\int_{[0,t]}\int_{\R\setminus\{0\}}g(x)(\mu^X(ds,dx)-\nu^X(ds,dx))\right)_{t\ge 0} \quad \mbox{is a local martingale}.
\]
In particular for any such function we get the compensation formula
\[
\E\Big( \sum_{t\le T} g(\Delta X_t) \Big)= \int g(x)\nu^X([0,T]\times dx)
\]
which extends to any nonnegative function $g$ or satisfying $\displaystyle \int_{\R\setminus\{0\}} |g(x)|\nu^X([0,T]\times dx)<+\infty$. The L\'evy measure $\nu^X$ satisfies 
\[
\nu^X(\{0\}\times \R)=\nu^X(\{t\}\times dx)=0,\;\quad \int_{\R} (x^2\wedge 1) \nu^X([0,t]\times dx)<+\infty, \; t\!\in\R_+.
\]

In what follows, we make the convention that $\Delta \alpha_{\infty}=0$ for any c\`adl\`ag function $\alpha$ defined on $\R_+$.

\smallskip First, owing to Lemma~\ref{LemTech}$(a)$ and the result about the standard Brownian  in Section~\ref{Trois}, we can assume that $c \equiv 0$ in what follows, $i.e.$ that  there is no Brownian component. Then we define  two independent marked Poisson processes with positive jumps as follows
\[
\widetilde X^{\varepsilon,\pm }_t = \sum_{s\le t }( \Delta X_s)_\pm \mbox{\bf 1}_{\{(\Delta X_s)_\pm> \varepsilon\}}
\] 
and $\widetilde X^{\varepsilon} = \widetilde X^{\varepsilon,+}-\widetilde X^{\varepsilon,- }$. For each process, we define their  inter jump sequence $(\widetilde \Theta^{\varepsilon,\pm}_n)_{n\ge 0}$  $i.e.$, with the convention $\widetilde \Theta^{\varepsilon,\pm}_0=0$, 
\[
\widetilde \Theta^{\varepsilon,\pm}_{n+1}=\min \{s>\widetilde S^{\pm}_n \,|\,    (\Delta X_{\widetilde S^{\varepsilon,\pm}_n+s})_\pm >\varepsilon      \}\!\in (0,+\infty],\quad n\ge 0,
\]
where $\widetilde S^{\varepsilon,\pm}_n = \widetilde \Theta^{\varepsilon,\pm}_1+\cdots+\widetilde \Theta^{\varepsilon,\pm}_n$. 

Both  processes $\widetilde X^{\varepsilon,+}$  and $\widetilde X^{\varepsilon,-}$ are independent since they have no common jumps. Furthermore the four sequences    $(\widetilde \Theta^{\varepsilon,\pm}_{n})_{n\ge 1}$ and $(\Delta \widetilde X^{\varepsilon,\pm}_{\widetilde S^{\pm}_n})_{n\ge 1}$ are mutually independent  and made of mutually independent terms. 

\smallskip
Let $F$ be a bounded measurable non-decreasing defined on functional on $I\!\!D([0,T],\R)$. Now, for every $n\ge   1$, we define   on $\R^{2n}$   the function $F_n$ by
\[
F_{n}(\xi_1,\theta_1, \dots,\xi_n,\theta_n)= F\Big(\big(\sum_{k=1}^n(\xi_k)_{_+} \,\mbox{\bf 1}_{\{\theta_1+\cdots+\theta_n\le t\}}\big)_{t\in [0,T]}\Big), \;  \xi_1,\ldots,\xi_n\!\in \R,\;\theta_1,\ldots,\theta_n\!\in \overline \R. 
\]
It is straightforward that the functions  $F_{n}$ are non-decreasing in each variable $\xi_i$ and non-increasing in each variable $\theta_i\!\in \overline \R_+$.

For every $n\ge 1$, set $\displaystyle \widetilde X^{\varepsilon,n,\pm}_t= \sum_{k=1}^n \big(\Delta X_{\widetilde S^{\varepsilon,\pm}_k}\big)_{\pm} \mbox{\bf 1}_{\{\widetilde S^{\varepsilon,\pm}_k\le t\}}$ so that 
 \[
F(\widetilde X^{\varepsilon,n,\pm})=  F_{n}\big(((\Delta X_{\widetilde S^{\varepsilon,\pm}_k})_{\pm} , \widetilde \Theta^{\varepsilon,\pm}_{k})_{k=1,\ldots,n}\big).
\]
Consequently, if $F$ and $G$ are co-monotone (measurable) functionals on  $I\!\!D([0,T],\R)$, it follows from Proposition~\ref{ProFinidim}$(b)$ (co-monotony principle for mutually independent random variables) that  
\[
\E\, F(\widetilde X^{\varepsilon,n,\pm})G(\widetilde X^{\varepsilon,n,\pm})\ge \E\, F(\widetilde X^{\varepsilon,n,\pm})\E\, G(\widetilde X^{\varepsilon,n,\pm}).
\]
Now 
\[
\sup_{t\in [0,T]}|\widetilde X^{\varepsilon,\pm}_t-\widetilde X^{\varepsilon,n,\pm}_t|\le \sum_{k\ge n+1} \big(\Delta X_{\widetilde S^{\varepsilon,\pm}_k}\big)_{\pm}\mbox{\bf 1}_{\{\widetilde S^{\varepsilon,\pm}_k\le T\}}
\]
so that 
\[
\P\big(\sup_{t\in [0,T]}|\widetilde X^{\varepsilon,\pm}_t-\widetilde X^{\varepsilon,n,\pm}_t|>0\big) \le \P(\widetilde S^{\pm}_{n+1} \le T)\to 0\quad \mbox{as } \quad n\to \infty
\]
since the process $X$  has finitely many jumps of size greater than $\varepsilon$ on any bounded time interval. The continuity of $F$ and $G$ transfers the co-monotony inequality  to $\widetilde X^{\varepsilon,\pm}$. In turn, the independence of these two processes, combined with Lemma~\ref{LemTech}$(a)$, propagates   co-monotony to the global Poisson process $\widetilde X^{\varepsilon}$.

\smallskip Noting that $F$ and $F(b^{\varepsilon}+.)$ have the same monotony (if any), one derives that $X-M^{\varepsilon}$ satisfies the co-monotony principle for every $\varepsilon>0$. One concludes by noting that $\|M^{\varepsilon}\|_{\sup}\to 0$ as $\varepsilon\to 0$ in $L^2$.~$\cqfd$

\subsection{C\`adl\`ag Markov processes and $\|\,.\,\|_{L_{_T}^p(\mu) }$-continuous functionals}\label{CadlagMarkov}
It is often convenient  to consider some path spaces of the form $L^p([0,T], \mu)$ where $\mu$ is a $\sigma$-finite measure and $p\!\in[1,+\infty)$, especially because of the properties of differentiation on these spaces which allow the natural introduction of gradient fields. Of course, less functionals are continuous for such a topology than with the $\|\,.\,\|_{\sup}$-norm topology when the process $X$ has continuous (or even c\`adl\`ag) paths.

Then, following the lines of the proof of Proposition~\ref{ContMarkov} but with a new canonical approximation procedure of a function $\alpha$, this time  by a stepwise constant function, one  shows the following property.
\begin{Pro}\label{ContMarkov} Let $(X_t)_{t\in [0,T]}$ be a c\`adl\`ag Markov process defined on a probability space $(\Omega,{\cal A},\P)$ with  transitions operators $(P_{s,t})_{t\ge s\ge 0}$ satisfying the monotony property. Let $\mu$ be a finite measure on $([0,T], {\cal B}or([0,T]))$ and let $p\!\in [1,+\infty)$. Let $F,G:\D([0,T], \R)\to \R$ be two $\mu$-co-monotone functionals, $\P_{_X}$-$a.s.$  continuous  with respect to the $L_{_T}^p(\mu) $-norm  on $\D([0,T], \R)$. If $F(X)$, $G(X)$ and $F(X)G(X)$ are  integrable or have $\P_{_X}$-$a.s.$ a common constant  sign, then
\[
\E\,F(X)G(X)\ge \E \, F(X) \E\, G(X).
\]
\end{Pro}

\ni {\bf Proof.} For every $\alpha\!\in \D([0,T], \R)$ and very integer $n\ge1$ we define  the stepwise constant approximation operator $\widetilde \alpha^{(n)}$ defined by~(\ref{stepwise}).  
It is clear that $\alpha^{(n)}(t)\to \alpha(t)$ at every   $t\!\in [0,T]$ and that the sequence $(\alpha^{(n)})_{n\ge 1}$ is bounded by $\|\alpha\|_{\sup}$. Hence $\alpha^{(n)}$ converges to $\alpha$ in every $L_{_T}^p(\mu) $, $1\le p<+\infty$. The rest of the proof is similar to that of Proposition~\ref{ContMarkov}.~$\cqfd$

\section{Applications}\label{Cinq}

\subsection{Functional antithetic simulation method}

Of course, the first natural application is a functional version of  the antithetic simulation method briefly described in the introduction.  To be more precise, let $X$ be a process taking values in a  vector subspace $E\subset {\cal F}([0,T], \R)$ (partially ordered by the pointwise order) satisfying a functional co-monotony principle in the sense of Definition~\ref{Def:2}. If, furthermore, it is invariant in distribution under a continuous  non-increasing mapping $T:E\to E$ (by non-increasing we mean that $\alpha\le \beta \Rightarrow T(\alpha)\ge T(\beta)$, $\alpha$, $\beta\!\in E$) then for any $\P_{_X}$-$a.s.$ sup-norm continuous  monotone  functional $F:E\to \R$ (square integrable or with constant sign)  
\[
{\rm Cov}( F(X), F\!\circ\! T(X)) = \E \big(  F(X)F\!\circ\! T(X)\big)-\big(\E\, F(X)\big)^2\le 0.
\] 
As a consequence, in order to compute $\E\,F(X)$ by a Monte Carlo simulation, it follows that  the computation of (independent copies of) $\displaystyle \frac {F(X)+F \!\circ\! T(X)}{2}$ will induce, for a prescribed simulation budget, a lower variance than a simulation only computing (independent copies of)   $F(X)$ like in the scalar framework. In practice such simulations rely on discretization schemes of $X$ for which the co-monotony principle is only true asymptotically (when the discretization step will go to zero). So is the case for the Euler scheme of a Brownian diffusion with non deterministic diffusion coefficient. 

It remains that this strongly suggests,  in order to compute $\E\,F(X)$ where $X=(X_t)_{t\in [0,T]}$ is  a Brownian diffusion (say the unique strong solution to an $SDE$ starting at $x$) and $F$ is a   sup-norm continuous  monotone  functional,   to simulate systematically two coupled  paths of an   Euler scheme (with a small enough step):  one with a sequence of Brownian increments $(W_{t^n_{k+1}}-W_{t^n_k})_{k\ge 0}$ and one with its opposite $-(W_{t^n_{k+1}}-W_{t^n_k})_{k\ge 0}$. In fact we know, $e.g.$ from~\cite{REYO} (chap.IX, p.341) that    $X=\Xi(W)$. Although we do not know whether $F\!\circ \! \Xi$ is monotone (and $\P_{_W}$-$a.s.$ sup-norm continuous)  the sign of the covariance can be roughly tested on a small simulated sample.

\subsection{A first application  to  peacocks}
The aim of this section is to prove that the (centered) antiderivative of an integrable process satisfying a co-monotony principle is a peacock in the sense of the definition given in the introduction.

\begin{Pro}\label{antider} Let $X=(X_t)_{t\ge 0}$ be an integrable    c\`adl\`ag process satisfying a co-monotony principle for the sup norm on every interval $[0,T]$, $T>0$, and let $\mu$ be a Borel measure on $(\R_+, {\cal B}or(\R_+))$.  Assume that $\sup_{[0,t]}\E\, |X_s|<+\infty$ for every $t>0$ and that $t\mapsto \E\, X_t $ is c\`adl\`ag. Then the process 
\[
\Big(\int_{[0,t]}(X_s-\E\, X_s)\mu(ds)\Big)_{t\ge 0} \mbox{ is a peacock.}
\]
\end{Pro}

\noindent {\bf Remark.} If $\displaystyle \sup_{t\in [0,T]} \E|X_t|^{1+\varepsilon} <+\infty$ for an $\varepsilon>0$, then $t\mapsto \E\, X_t $ is c\`adl\`ag by a uniform integrability argument. 

\bigskip
\noindent {\bf Proof.} First we may assume without loss of generality that the process $X$ is centered since $(X_t-\E\,X_t)_{t\in [0,T]}$ clearly satisfies a co-monotony principles on $I\!\!D([0,T], \R)$ for every $T>0$   if $X$ does (since $t\mapsto \E\, X_t $ is c\`adla\`ag). Set for convenience $Y_t = \int_0^t X_s \mu(ds)$. It is clear from the assumption and Fubini's Theorem that $Y_t\!\in L^1(\P)$ and $\E\,Y_t=0$. 

\smallskip
\noindent   $\rhd$ {\sc Step~1~:} Let $\varphi:\R\to \R$ be a convex function with linear growth (so that $\varphi(Y_t)\!\in L^1(\P)$ for every $t\ge 0$). Its right derivative $\varphi'_r$ is a non-decreasing bounded function. The convexity of the function  $\varphi$ implies, for every $x$, $y\!\in \R$, 
\[
\varphi(y)-\varphi(x)\ge \varphi'_r(x)(y-x)
\]
so that, if $t_1<t_2$
\[
\varphi(Y_{t_2})-\varphi(Y_{t_1})\ge \Phi_{t_1}(X)\int_{(t_1,t_2]} X_s\,\mu(ds)
\]
where $\Phi_{t_1}(\alpha)=   \varphi'_{r}\Big(\int_{[0,t_1]}\alpha(s)\mu(ds)\Big)$ is bounded, continuous for the sup norm topology on $I\!\!D([0,t_2], \R)$ and pointwise non-decreasing  for the pointwise order. The functional $\alpha\mapsto \int_{(t_1,t_2]}\alpha(s)\mu(ds)$ is also continuous for the  sup norm topology, pointwise non-decreasing and 
$$
\Big|\E \int_{(t_1, t_2]}X_s\mu(ds)\Big| \le \E \int_{(t_1, t_2]}|X_s|\,\mu(ds)\le  \mu((t_1,t_2])\sup_{s\in [t_1,t_2]} \E|X_s|<+\infty.
$$
 Consequently, owing to the co-monotony principle, we get 
\[
\E \Big(\Phi_{t_1}(X)  \int_{(t_1,t_2]} \hskip -0.5 cm X_s \mu(ds)\Big) \ge \E\,  \Phi_{t_1}(X) \E\Big(\int_{(t_1,t_2]}\hskip -0.5 cm  X_s \mu(ds) \Big)=  \E\,  \Phi_{_{T}}(X) \times \int_{(t_1,t_2]}\hskip -0.5 cm  \E\, X_s \mu(ds)  =0
\]
so that $\E\, \varphi(Y_{t_2}) \ge \E\,\varphi(Y_{t_1})\!\in L^1(\P)$.

\medskip
\noindent {\sc Step~2}: Assume now that $\varphi$ is simply convex. 
For every $A>0$, we define the following convex function $\varphi_{_A}$ with linear growth:
\[
\varphi_{_A}(y) = \left\{\begin{array}{ll}\varphi(y) & \mbox{if } |y|\le A\\ \varphi(A) +\varphi'_r(A)(y-A)& \mbox{if } y\ge A\\
\varphi(-A) +\varphi'_r(-A)(y+A) &  \mbox{if } y\ge A.\end{array} \right.
\]

It is clear that  $\varphi_{_A}\uparrow \varphi$ as $A\uparrow +\infty$ and that $\E\,\varphi_{_A}(Y_{t_2}) \ge \E\,\varphi_{_A}(Y_{t_1})$ by Step~1.


Now $\varphi_{_A}$ has linear growth so that $\varphi_{_A}(Y_t)\!\in L^1(\P)$. Consequently it follows from the monotone convergence theorem that $\E\, \varphi_{_A}(Y_t) \uparrow \E\, \varphi(Y_t)\!\in (-\infty,+\infty]$ as $A\uparrow +\infty$. This completes the proof.$\cqfd$

\subsection{From the sensitivity of Asian path-dependent options to peacocks} 

 Let $\mu$ be a finite measure on $([0,T], {\cal B}or([0,T]))$ and, for every $p\!\in [1,+\infty)$, let $q$ denote its H\"older conjugate.  Note that, of course,   $\D([0,T], \R)\subset L^{\infty}_{_T}(\mu) \subset \cap_{p\ge 1} L^p_{_T}(\mu)$. 
 
 \begin{Dfn} $(a)$ Let $p\!\in [1,+\infty)$. A  measurable functional  $F:L_{_T}^p(\mu) \to \R$ is {\em regularly differentiable} on $\D([0,T], \R)$ if, for every  $\alpha\!\in \D([0,T], \R)$, there exists a measurable ``gradient" functional $\nabla F:\big([0,T]\times \D([0,T], \R), {\cal B}or([0,T])\otimes {\cal D}_T\big)\to \R$ such that
 \begin{equation}\label{DiffL^p}
\left\{\begin{array}{ll}(i)& \nabla F(.,\alpha)\!\in L_{_T}^q(\mu)\\
(ii) &\lim_{\|h\|_{L_{_T}^p(\mu) }\to 0, h \in L_{_T}^p(\mu)}\frac{\Big|F(\alpha+h)-F(\alpha) -\int_0^T \nabla F(s,\alpha)h(s)\mu(ds)\Big| }{\|h\|_{L_{_T}^p(\mu) }}=0.
\end{array}\right.
 \end{equation}
\noindent $(b)$ Furthermore, a gradient functional  $\nabla F$ is {\em monotone} if, for every $t\!\in [0,T]$, $\nabla F(t,.)$ is monotone on $\D([0,T], \R)$ and if  this monotony does not depend on $t\!\in [ 0,T]$.
\end{Dfn}

\begin{Pro}\label{prepeacock}Let $X=(X_t)_{t\in [0,T]}$ be a (c\`adl\`ag) {\em PII} such that, for every $u\!\in \R$, $L(u,t)=\E \,e^{uX_t}$ is bounded and bounded away from $0$ over $[0,T]$ so that, in particular,  the function $\Psi(u,t)= \log \E\, e^{uX_t}$ can be defined as a real valued function. Let  $F:L_{_T}^p(\mu) \to \R$ be a measurable functional, regularly differentiable  with a monotone gradient $\nabla F$ on $ \D([0,T], \R)$  satisfying the following Lipschitz continuity assumption
\[
\forall\, \alpha,\, \beta \!\in \D([0,T], \R),\quad | F(\alpha)- F(\beta)| \le[F]_{\rm Lip} \|\alpha-\beta\|_{L_{_T}^p(\mu) }.
\]
Set, for every $\sigma>0$,  
 \begin{equation}\label{def:f}
 f(\sigma)=\E \left( F\Big(e^{\sigma X_.-\Psi(\sigma,.)}\Big)\right) .
 \end{equation}

Then, under the above assumptions, the function $f$ is (differentiable and) non-decreasing. 
\end{Pro}

\ni {\bf Remark.} At least for L\'evy processes, the assumption $\sup_{t\in [0,T]}\E e^{uX_t}<+\infty$, $u\!\in \R$, is satisfied as soon as $\E \,e^{uX_t}<+\infty$ for every $u\!\in \R$ (see~\cite{SATO}, Theorem 25.18, p.166).

\bigskip
Before proving the proposition, we need the following technical lemma about the regularity of function $L$ whose details of proof are left to the reader.

\begin{Lem}\label{TechL} Under the assumption made on the function $L$ in Proposition~\ref{prepeacock}, the function $\Lambda$ defined on $\R_+^2$ by  $\Lambda(a,t) = \E \,e^{a|X_t|}$ is finite. Then for every $a\!\in (0,+\infty)$, $L$ is Lipschitz continuous in $u$ on $[-a,a]$, uniformly in $t\!\in [0,T]$, with Lipschitz coefficient (upper-bounded by) $\Lambda(a,T)$. Furthermore, for every $u\!\in \R$,  there exists $\kappa_{u,T}>0$ and $\varepsilon=\varepsilon(u,T)>0$ such that
\[
\forall\, t\!\in [0,T], \; \forall\, u'\!\in [u-\varepsilon, u+\varepsilon],\; L(u,t)\ge \kappa_{u,T}.
\]
\end{Lem}

\ni {\bf Proof of Proposition~\ref{prepeacock}.} Formally, the derivative of $f$ reads
\begin{eqnarray*}
f'(\sigma)&=& \E\left(\int_0^T \nabla F\Big(e^{\sigma X_.-\Psi(\sigma,.)},t\Big)e^{\sigma X_t-\Psi(\sigma,t)}(X_t-\Psi'_{\sigma}(\sigma,t)\big)\mu(dt) \right)\\
&=& \int_0^T \E\left(  \nabla F\Big(e^{\sigma X_.-\Psi(\sigma,.)},t\Big)e^{\sigma X_t-\Psi(\sigma,t)}(X_t-\Psi'_{\sigma}(\sigma,t)\big)\right)\mu(dt).
\end{eqnarray*}
To justify that  this interchange of differentiation and expectation in the first line  is valid we need to prove that the ratio \[
\frac{ F\Big(e^{\sigma' X_.-\Psi(\sigma',.)}\Big)- F\Big(e^{\sigma X_.-\Psi(\sigma,.)}\Big)}{\sigma'-\sigma},\;\sigma'\neq \sigma,\;\sigma,\, \sigma'\!\in[\epsilon_0,1/\epsilon_0], \;\epsilon_0>0,
\]
 is $L^{1+\eta}$-bounded for an $\eta>0$. Without loss of generality, we may assume that $p=1+\eta>1$ since $\|\,.\,\|_{L^p_{_T}}\le \mu([0,T])^{\frac 1p-\frac{1}{p'}}\|\,.\,\|_{L^{p'}_{_T}}$ if $1\le p\le p'$. This follows from the Lipschitz continuity of $F$ and from  the properties of the Laplace transform $L$ established in Lemma~\ref{TechL}.


\smallskip Let ${\cal G}_t := \sigma(X_s-X_t, \, s\!\in [t,T])$. This $\sigma$-field  is independent of ${\cal F}^X_t$.  Elementary computations show that, for every $t\!\in [0,T]$, 
\[
\E\left(  \nabla F\Big(e^{\sigma X_.-\Psi(\sigma,.)},t\Big) \,|\, {\cal G}_t\right)=\Phi\Big(X_t\!-\!\frac{\Psi(\sigma,.)}{\sigma},\big(X_s\!-\!\frac{\Psi(\sigma,.)}{\sigma}\big)_{s\in [0,t]},t\Big)
\]
where, for every $\beta\!\in I\!\!D([0,t], \R)$,
\[
 \Phi(\xi,\beta,t) = \E\left(\nabla F\big(e^{\sigma(X_.-X_t)-(\Psi(\sigma,.)-\Psi(\sigma,t))+\sigma \beta(t)}\mbox{\bf 1}_{(t,T]}+ e^{\sigma \beta}\mbox{\bf 1}_{[0,t]},t\big)\right).
\]

Note  that, for every $t\!\in[0,T]$, the function $\Phi(.,.,t)$ is non-decreasing in both remaining arguments. Now
\begin{equation}\label{eq:Repres}
f'(\sigma) = \int_0^T\hskip-0.25 cm  \E\! \left( \!\Phi\Big(\big(X_s\!-\!\frac{\Psi(\sigma,.)}{\sigma}\big)_{s\in [0,t]},t\Big)
e^{\sigma (X_t-\frac{\Psi(\sigma,t)}{\sigma})}\big(X_t-\Psi'_{\sigma}(\sigma,t)\big)\!\right)\!\mu(dt). 
\end{equation}

Set $\Q_{(t)} = e^{\sigma X_t-\Psi(\sigma,t)}.\P$. It is classical that $(X_s)_{s\in [0,t]}$ is still a {\em PII} under $\Q_{(t)}$ with exponential moment at any order and a $\log$-Laplace transform $\Psi_{(t)}$ given by 
$$
\Psi_{(t)}(u,s) = \Psi(\sigma+u,s)-\Psi(\sigma,s).
$$ 
Note that  $\Psi_{(t)} $ does not depend on $t$ but   on $\sigma$. Consequently, for every $s\!\in [0,t]$, 
\[
\E_{\Q_{(t)}}\big(X_s\big)= \frac{\partial \Psi_{(t)}}{\partial u}(0,s)= \Psi'_{\sigma}(\sigma,s)
\]
where $\Psi'_{\sigma}(\sigma,s)$ denotes the partial derivative of $\Psi$ with respect to $\sigma$. Putting $\widetilde X_s= X_s-\Psi'_{\sigma}(\sigma,s)$, we get

\begin{equation}\label{eq:Repres}
f'(\sigma) = \int_0^T \E_{\Q_{(t)}} \left(  \Phi\Big(\big(\widetilde X_s + \Psi'_{\sigma}(\sigma,s)-\frac{\Psi(\sigma,s)}{\sigma}\big)_{s\in [0,t]},t\Big)\widetilde X_t\right)\mu(dt).
\end{equation}

Applying the co-monotony principle to the {\em PII} $\widetilde X$ and to the two non-decreasing $L_{_T}^p(\mu) $-continuous functionals $F(\alpha)= \Phi\big(\big(\alpha(s) + \Psi'_{\sigma}(\sigma,s)-\frac{\Psi(\sigma,s)}{\sigma}\big)_{s\in [0,t]}\big)$ and $G(\alpha)=\alpha(t)$ yields that, for every $t\!\in [0,T]$,  
\[
 \E_{\Q_{(t)}} \left(  \Phi\Big( \big(\widetilde X_s + \Psi'_{\sigma}(\sigma,s)-\frac{\Psi(\sigma,s)}{\sigma}\big)_{s\in [0,t]},t\Big) \widetilde X_t\right)\ge 0
\]
since $ \E_{\Q_{(t)}} \, \widetilde X_t=0$. As a consequence, $f$ is a non-decreasing function.~$\cqfd$

\begin{Cor}\label{Cor:peacock}  Under the assumptions of Proposition~\ref{prepeacock} on the c\`adl\`ag {\em PII} $X$, the process \\ $\displaystyle  \Big( \int_0^Te^{\sigma X_t-\Psi(\sigma,t)}\mu(dt)\Big)_{\sigma\!\in \R_+}$  is a peacock.
\end{Cor} 

\medskip
\ni{\bf Proof.}  Let $\varphi:\R\to\R$ be   a convex function and, for every $A>0$, let $\varphi_A$ be defined by $\varphi_A'(x)=\varphi(x) $ if $x\!\in [-A,A]$ and $\varphi_A$ affine and differentiable on $(-\infty,-A]\cup [A,+\infty)$.  It is clear that $\varphi_A \uparrow \varphi$ since $\varphi $ takes values in $(-\infty, +\infty]$. Then set $\varphi_{A,\varepsilon}(x) = \E\, \varphi_A(x+\varepsilon Z)$ where $Z\stackrel{\cal  L}{\sim}{\cal N}(0;1)$. The function $\varphi_{A,\varepsilon}$ is (finite) convex, infinitely differentiable,  Lipschitz continuous and  converges uniformly to $\varphi_A$ when $\varepsilon\to 0$. 
The functional  $F_{A,\varepsilon}(\alpha)= \varphi_{A,\varepsilon}\big(\int_0^T\alpha(t)\mu(dt)\big)$ satisfies the assumptions of the above Proposition~\ref{prepeacock} so that the function $f_{A,\varepsilon} $ defined by~(\ref{def:f})  is non-decreasing. Letting $\varepsilon\to 0$ and   $A\to +\infty$ successively implies that the function $f$ related (still through~(\ref{def:f})) to the original functional $F(\alpha)= \varphi\big(\int_0^T\alpha(t)\mu(dt)\big)$      is non-decreasing which completes the proof (for $A\uparrow +\infty$ the arguments are those of  the proof of Proposition~\ref{antider}.~$\cqfd$

\bigskip
\ni {\bf Remarks.} $\bullet$ In fact this proof remains close in spirit to that proposed in~\cite{YORetal}. Roughly speaking we replace the notion of {\em conditional monotony} used in~\cite{YORetal} by a functional co-monotony argument (which also spares a time discretization phase). The notion of conditional monotony  and its applications have been extensively investigated  in the  recent PhD thesis of A.  Bogso (see~\cite{BOG}).  Conditional monotony has been developed on the basis of finite dimensional distributions of a process but it is clear that a functional version can be derived for (continuous) functionals. Then, when the parameter of interest is time,  the connection with functional co-monotony looks clear since it corresponds to a {\em weak form}  of the functional co-monotony principle restricted to couples of functionals of the form 
$F(\alpha^t)$ and $G(\alpha)= g(\alpha(t))$ ($\alpha^t$ denotes the stopped function $\alpha$ at $t$).

\medskip
\ni $\bullet$ As already noticed in~\cite{YORetal}, specifying $\mu$ into $\delta_T$ or $\frac 1T \lambda_{|[0,T]}$ provides  the two main results for peacocks devised from $e^{\sigma X_t-\Psi(\sigma,t)}$. When $\mu =\frac 1T \lambda_{|[0,T]}$ one can combine the above results with some self-similarity property of the {\em PII} process  $(X_t)_{t\in [0,T]}$ (if any)  to produce other peacocks. So is the case with the seminal example investigated in~\cite{Carretal} where the original aim was, for financial purposes,  to prove that
\[
\left(\frac 1t \int_0^t e^{B_s-\frac s2}ds\right)_{t \in(0,T]}\mbox{ is a peacock}.
\]

Many other examples of this type are detailed in~\cite{YORetal, BOG}.

\bigskip
\noindent {\sc Application to a class of Asian options.} As concerns the sensitivity of exotic derivatives, one can derive or retrieve classical results in a Black-Scholes model for the class of Asian options with convex payoff. To be precise, we consider  payoff functionals  of the form $\Phi_T= \varphi \big(\frac 1T\int_0^TS_sds\big)$ where $\varphi$ is a nonnegative convex function (with linear growth) and $S_t= s_0 e^{(r-\frac{\sigma^2}{2})t+\sigma W_t}$, $t\!\in [0,T]$,  where $s_0>0$, $\sigma>0$ and $W$ is a standard Brownian motion ($r$ is a possibly negative interest rate). The holder of an option contract ``written" on this payoff receives in cash at the {\em maturity} $T>0$ the value of the payoff $\Phi_{_T}$. Classical arbitrage arguments yield that the premium or price at time $0$ of such an option is given by 
\[
\mbox{Premium}_0 (s_0,\sigma, r,T) = e^{-rT} \E \,\varphi \Big(\frac 1T\int_0^TS_sds\Big).
\]

By considering the measure $\mu(dt)= e^{rt}\frac 1T \lambda_{|[0,T]}(dt)$, one derives from Corollary~\ref{Cor:peacock} that $\sigma\mapsto \mbox{Premium}_0 (s_0,\sigma, r,T) $ is non-decreasing. When $r=0$ a change of  variable shows that the premium is also  non-decreasing as a function of the  maturity $T$.

\subsection{Application to quasi-universal bounds for  barrier options}\label{bounds}

Let $S=(S_t)_{t\in [0,T]}$ be a c\`adl\`ag nonnegative stochastic process  defined on a probability space $(\Omega,{\cal A}, \P)$, modeling the price dynamics of a risky asset. We will assume that $\P$ is a pricing measure in the sense that  derivatives products ``written" on the asset $S$ are priced under $\P$. In particular we do not ask $\P$ to be risk-neutral. We assume for convenience that the zero-coupon bond (also known as the riskless asset) is constant equal to $1$ (or equivalently that all interest rates are constant equal to $0$) but what follows remains true if  the price dynamics of this bond   is deterministic.

 \smallskip For notational convenience, for a c\`adl\`ag function $\alpha:[0,T]\to \R$, we   will denote by $_*\alpha_t :=\inf_{s\le t} \alpha_s$, $(t\!\in [0,T])$,  the running minimum of the function $\alpha$ and by  $\alpha^*_t :=\sup_{s\le t} \alpha_s$ its running maximum process. In what follows we will extensively use the following classical facts: $\alpha\mapsto _*\!\!\alpha_{_T}$ and $\alpha\mapsto \alpha^*_{_T}$ are Skorokhod continuous on $I\!\!D([0,T], \R)$ and, for every $t\!\in [0,T)$,  $\alpha\mapsto _*\!\alpha_t$ and $\alpha\mapsto \alpha^*_t$ are sup-norm continuous  on $I\!\!D([0,T], \R)$ (and Skorokhod continuous at every $\alpha$ continuous at $t$).

 We assume  throughout this section that {\em the asset price dynamics $(S_t)_{t\in [0,T]}$ satisfies  a functional co-monotony principle}. This seems is a quite natural and general assumption given the various classes of examples detailed above.


\medskip 
\ni  We will focus on   {\em Down-and-In Call} and {\em Down-and-Out Call} with maturity $T$. The payoff functional of a  {\em Down-and-In Call} with maturity $T$  is defined for every strike price $K>0$ and every barrier $L\!\in (0,S_0)$ by
\[
F_{D\&I}(\alpha)= \big(\alpha(T)-K)_+\mbox{\bf 1}_{\{_*\!\alpha_{_T}\le L\}}.
\]
This means that the holder  of the    {\em Down-and-In Call}  written on the risky asset $S$   contract receives $F_{D\&I}(S)$ at the maturity $T$, namely $S_T-K$ at the maturity $T>0$ provided this flow is positive and  the asset attained at least once the (low) barrier $L$ between $0$ and $T>0$.

The premium  of such a contract at time $0$ is defined by  
$$
{\rm Call}_{D\&In}(K,L,T) = \E \big(F_{D\&I}(S)  \big).
$$
  We will denote by $ {\rm Call}(K,T) = \E\, (S_{_T}-K)_+$ the  premium of the regular (or vanilla) {\em Call} option with strike $K$ (and maturity $T$). 

\begin{Pro}If  the nonnegative c\`adl\`ag process $(S_t)_{t\in [0,T]}$ satisfies  a finite dimensional co-monotony principle (hence  functional co-monotony principle on $I\!\!D([0,T], \R)$ for the Skorokhod topology), then the following semi-universal bound holds:
\[
 {\rm Call}_{D\&In}(K,L,T)\le {\rm Call}(K,T)  \P(_*S_T\le L).
\]
\end{Pro}

\noindent {\bf Proof.}    For every $\alpha\!\in I\!\!D([0,T], \R)$ and every $\varepsilon>0$, we have
\[
F_{D\&I}(\alpha)\le  \big(\alpha(T)-K)_+\left (\Big(1-\frac{_*\alpha_T -L}{\varepsilon}\Big)_+\wedge 1\right).
\]
The two functionals involved in the product of the right hand side of the above equation are clearly anti-monotone, nonnegative and continuous with respect to the Skorokhod topology, consequently
\[
 {\rm Call}_{D\&In}(K,L,T) \le  {\rm Call}(K,T)\E\left(\Big(1-\frac{_*S_T -L}{\varepsilon}\Big)_+\wedge 1\right).
\]
The result follows by letting $\varepsilon \to 0$ owing to Fatou's Lemma.$\cqfd$

\bigskip
As concerns the  {\em Down-and-Out Call} with payoff
\[
F_{D\&O}(\alpha)= \big(\alpha(T)-K)_+\mbox{\bf 1}_{\{_*\alpha_T >L\}}
\]
for which the holder of the option receives $S_T-K$ at the maturity $T>0$  if this flow is positive and if the asset did not attain the (low) level $L$ between $0$ and $T>0$, one gets, either by a direct approach or by 
using  the obvious parity equation ${\rm Call}_{D\&In}(K,L,T)+{\rm Call}_{D\&Out}(K,L,T)={\rm Call}(K,T)$, 
\[
  {\rm Call}_{D\&Out}(K,H,T)\ge {\rm Call}(K,T)\P(_*S_T> L).
\]

 Similar bounds can be derived for   {\em Up-and-In} and {\em Up-and-Out Calls} with barrier $L>S_0$ (and strike $K$), namely
 \[
  {\rm Call}_{U\&In}(K,L,T)\ge {\rm Call}(K,T)\P(S^*_T>  L )
 \] 
 and
  \[
   {\rm Call}_{U\&Out}(K,L,T)\le {\rm Call}(K,T)\P(S^*_T\le L ).
 \]

If one considers extensions of the above payoffs in which the barrier needs to be (un-)knocked strictly prior to $T$, at a time $T'<T$, similar semi-universal bounds can be obtained provided one of the following assumption is true $\P(S_{T'}=S_{T'-})=1$ or $(S_t)_{t\in [0,T]}$ satisfies a functional co-monotony principle with respect to the sup-norm on $I\!\!D([0,T], \R)$.

\subsection{A  remark on running extrema}
If a  c\`adl\`ag  process $X=(X_t)_{t\in [0,T]}$ satisfies a co-monotony principle and $X_{_T}$ and $ \sup_{t\in[0,T]}X_t$ have no atom so that, for every $x$, $y\!\in \R$, $x\le y$, $\alpha\mapsto \big(\mbox{\bf 1}_{\{\alpha(T)\ge x\}},\mbox{\bf 1}_{\{\sup_{t\in [0,T]}\alpha(t)\ge y\}}\big)$ is $\P_{_X}$-$a.s.$ $\|\,.\,\|_{\sup}$-continuous,  then
\[
\forall\, y\!\in \R,\quad \P\big(\sup_{t\in [0,T]} X_t \ge y\big)= \inf _{x\le y} \P\big(\sup_{t\in [0,T]} X_t \ge y\,|\, X_{_T}\ge x\big).
\]

Of course the list of possible applications is not exhaustive. In more specific problems, one can take advantage of the functional co-monotony principle to establish less conservative inequalities and bounds on parameters of a problem.  A typical example is provided by~\cite{LALEPA} devoted to optimal order execution on a financial market, cited here since it  was partially at the origin of the present work.

\bigskip
\noindent {\bf Acknowledgement:} The author thanks the anonymous referee for  constructive suggestions and F.~Panloup for his help.


\begin{thebibliography}{}

\end{thebibliography}


\small  \begin{thebibliography}{0}
\bibitem{BAEM} {\sc P. Bauman, M. \'Emery} (2008). Peut-on ``voir" dans l'espace \`a $n$ dimensions, {\em  L'ouvert}, {\bf 116}:1-8.
\bibitem{BIL} {\sc P. Billingsley} (1999). {\em Convergence of probability measures}. Second edition. Wiley Series in Probability and Statistics: Probability and Statistics. A Wiley-Interscience Publication. John Wiley \& Sons, Inc., New York, 277p.
\bibitem{BOG}{\sc A.M. Bogso} (2012). {\em \'Etude de peacocks sous des hypoth\`eses de monotonie conditionnelle et de positivit\'e totale}, th\`ese de l'Universit\'e de Lorraine, 147p.
\bibitem{Carretal} {\sc P. Carr, C.-O.Ewald, Y. Xiao} (2008). On the qualitative effect of volatility and duration on prices of Asian options, {\em  Finance Research Letters}, {\bf 5}, 162-171.
\bibitem{DIA}{\sc P.H. Diananda} (1962). On nonnegative forms in real variables some or all of which
are nonnegative, {\em Proc. Cambridge Philos. Soc.},  {\bf 58}:17-25.
\bibitem{ENG} {\sc H. J. Engelbert }(1991).  On the theorem of T. Yamada and S. Watanabe, {\em Stochastics and Stochastic Reports}, {\bf 36}(3-4):205-216. DOI: 10.1080/17442509108833718
\bibitem{HAL}{\sc M. Hall} (1962).  Discrete problems, in {\em A Survey of Numerical Analysis}, John Todd
ed., McGraw-Hill.
\bibitem{HIRO} {\sc F. Hirsch, B. Roynette} (2012). A new proof of Kellerer Theorem,  {\em ESAIM : P\&S}, {\bf 16}:48-60.
\bibitem{YORetal} {\sc F. Hirsch, C. Profeta, B. Roynette, M. Yor} (2011).  {\em  Peacocks and Associated Martingales, With Explicit Constructions}.  Bocconi \& Springer, 430p.
\bibitem{JAC} {\sc J. Jacod} (2004). The Euler scheme for L\'evy driven stochastic differential equations: limit theorems, {\em Ann. Probab.}, {\bf  32}(3) (2004), 1830-1872.
\bibitem{JASH}{\sc J. Jacod, A.N. Shiryaev} (2003). {\em Limit theorems for stochastic processes}. Second edition. Grundlehren der Mathematischen Wissenschaften [Fundamental Principles of Mathematical Sciences],  {\bf 288}. Springer-Verlag, Berlin, 661p. 
\bibitem{JDPEPI}{\sc K. Joag-Dev, M. D. Perlman, L.D. Pitt} (1983).  Association of normal random variables and Slepian's inequality. {\em Ann. Probab.}, {\bf  11}(2):451-455. 
\bibitem{KASH} {\sc I. Karatzas,  S.E. Shreve} (1991). {\em  Brownian motion and stochastic calculus}, Second edition. Graduate Texts in Mathematics, 113. Springer-Verlag, New York,  470p.
\bibitem{KEL} {\sc   H.G. Kellerer} (1972). Markov-Komposition und eine Anwendung auf Martingale. (German) Math. Ann. {\bf 198}:99-122.
\bibitem{LALEPA} {\sc S. Laruelle,  C.-A. Lehalle, G. Pag\`es} (2011). Optimal posting price of limit orders: learning by trading,  submitted to {\em Mathematics \& Financial Economics}, pr\'epub-PMA-148, 2011.
\bibitem{LETA} {\sc M. Ledoux, M. Talagrand} (1991). {\em Probability in Banach spaces. Isoperimetry and processes}, Ergebnisse der Mathematik und ihrer Grenzgebiete (3) [Results in Mathematics and Related Areas (3)], {\bf  23}, Springer-Verlag, Berlin, 480p.
\bibitem{LUPA5} {\sc H. Luschgy, G. Pag\`es} (2009). Expansions for Gaussian processes and Parseval frames, {\em Electron. J. Probab.}, {\bf  14 }(42):1198-1221.
\bibitem{PIT} {\sc L.D. Pitt} (1982). Positively correlated normal variables are associated. {\em Ann. Probab.} {\bf 10}(2):496-499.
\bibitem{REYO} {\sc D. Revuz, M. Yor} (1999). {\em Continuous martingales and Brownian motion}, third edition. Grundlehren der Mathematischen Wissenschaften [Fundamental Principles of Mathematical Sciences], 293. Springer-Verlag, Berlin, 560p. 
\bibitem{SATO}{\sc K.I. Sato} (1999). {\em L\'evy Distributions and Infinitely Divisible Distributions}, Cambridge Studies in Advanced Mathematics,  Cambridge University Press, UK, 486p (first japanese edition in 1990).
\bibitem{ZBRA} {\sc  G. Zbaganu, M. Radulescu} (2009).
Trading prices when the initial wealth is random. {\em Proc. Rom. Acad. Ser. A Math. Phys. Tech. Sci. Inf. Sci.}, {\bf 10}(1):1-8.  
\end{thebibliography}
\end{document}